\documentclass[11pt,oneside,a4paper,leqno,usenames,dvipsnames]{amsart}
\usepackage{amsmath,amsfonts,amsthm}
\usepackage[utf8]{inputenc}
\usepackage[T1]{fontenc}
\usepackage{setspace}
\usepackage{enumitem}
\usepackage{amssymb}
\usepackage{mathtools}
\usepackage{hyperref}
\usepackage{appendix}
\usepackage{xcolor}
\usepackage[all]{xy}
\usepackage{vmargin}
\usepackage{tikz}
\usepackage{tikz-cd}
\usepackage{comment}
\usetikzlibrary{positioning}
\usepackage[document]{ragged2e}

\setmarginsrb{4.05cm}{3.5cm}{4.05cm}{3cm}{0pt}{1cm}{0pt}{1cm}
\newtheorem{theorem}{Theorem}[section]
\newtheorem{corollaire}[theorem]{Corollary}
\newtheorem{lemma}[theorem]{Lemma}
\newtheorem*{theorem*}{Theorem}
\newtheorem*{lemma*}{Lemma}
\newtheorem*{corollaire*}{Corollary}
\newtheorem*{question}{Question}
\newtheorem{prop}[theorem]{Proposition}

\theoremstyle{definition}
\newtheorem{definition}[theorem]{Definition}
\theoremstyle{remark}
\newtheorem{example}[theorem]{Example}
\newtheorem{remark}[theorem]{Remark}

\title[Coincidences of Division Fields]{Coincidences of Division Fields of an elliptic curve defined over a number field}
\author{Zoé Yvon}
\address{
Zoé Yvon, Institut de Mathématiques de Marseille (UMR 7373)
Site Sud, Campus de Luminy
Case 930
13288 MARSEILLE Cedex 9}
\email{zoenovy@free.fr}
\date{\today}

\newcommand\End{\operatorname{\mathrm{End}}}
\newcommand\Aut{\operatorname{\mathrm{Aut}}}
\newcommand\Gal{\operatorname{\mathrm{Gal}}}

\newcommand\got\mathfrak

\newcommand{\Z}{\mathbb{Z}}
\newcommand{\Q}{\mathbb{Q}}

\newcommand{\PP}{\mathbb{P}}

\newcommand{\K}{F}
\newcommand{\id}{\operatorname{\mathrm{id}}}

\newcommand{\D}{\operatorname{\mathrm{D}}}

\newcommand{\barr}{\overline}
\newcommand{\ord}{\mathrm{ord}}
\newcommand{\SL}{\mathrm{SL}_2}
\newcommand{\GL}{\mathrm{GL}_2}
\newcommand{\Mod}[1]{\ (\mathrm{mod}\ #1)}

\def\restriction#1#2{\mathchoice
              {\setbox1\hbox{${\displaystyle #1}_{\scriptstyle #2}$}
              \restrictionaux{#1}{#2}}
              {\setbox1\hbox{${\textstyle #1}_{\scriptstyle #2}$}
              \restrictionaux{#1}{#2}}
              {\setbox1\hbox{${\scriptstyle #1}_{\scriptscriptstyle #2}$}
              \restrictionaux{#1}{#2}}
              {\setbox1\hbox{${\scriptscriptstyle #1}_{\scriptscriptstyle #2}$}
              \restrictionaux{#1}{#2}}}
\def\restrictionaux#1#2{{#1\,\smash{\vrule height .8\ht1 depth .85\dp1}}_{\,#2}}
\usepackage{cancel}

\begin{document}

\justifying

\begin{abstract}
For an elliptic curve defined over a number field, the absolute Galois group acts on the group of torsion points of the elliptic curve, giving rise to a Galois representation in $\GL(\hat{\Z})$. The obstructions to the surjectivity of this representation are either local (\emph{i.e.} at a prime), or due to nonsurjectivity on the product of local Galois images. In this article, we study an extreme case: the coincidence \emph{i.e.} the equality of $n$-division fields, generated by the $n$-torsion points, attached to different positive integers~$n$. We give necessary conditions for coincidences, dealing separately with vertical coincidences, at a given prime, and horizontal coincidences, across multiple primes, in particular when the Galois group on the $n$-torsion contains the special linear group. We also give a non-trivial construction for coincidences not occurring over~$\Q$.
\end{abstract}

\maketitle


\section*{Introduction}


Let $F$ be a number field, $\overline{F}$ an algebraic closure of $F$ and $G_F=\Gal(\overline{F}/F)$ its absolute Galois group. Let $E/F$ be an elliptic curve. We know that the absolute Galois group $G_F$ acts on the group $E_\mathrm{tors}$ of the torsion points of $E/F$ encoded by the Galois representation \[\rho_E:G_F\to\Aut(E_\mathrm{tors})\simeq\GL(\hat{\Z}).\]
Serre \cite{Serre} proved that, if $E/F$ does not have CM, then the image of $\rho_E$ has finite index in $\GL(\hat{\Z})$. If $E/F$ has CM, then there exists a choice of basis for $E_\mathrm{tors}$ such that $\rho_E(G_F)$ is contained with finite index in a subgroup $N_{\delta,\phi}$, only depending on $\End(E)$ as order in a quadratic field, see \cite[Theorem 1.2]{lozrobcm}. Set $G=\GL(\hat{\Z})$ if $E/F$ does not have CM and $G=N_{\delta,\phi}$ if $E/F$ has CM. Let $G(p^\infty)$, respectively $G(m)$, be the image of $G$ in $\GL(\Z_p)$, respectively in $\GL(\Z/m\Z)$. Therefore, the mod $m$ Galois representation \[\rho_{E,m}:G_F\to G(m)\]is surjective for all $m$ coprime to a fixed integer. For a positive integer $m$, the non surjectivity of $\rho_{E,m}$ can be explained by two phenomena:\begin{enumerate}
    \item \textbf{Local conditions}: The non surjectivity of the $p$-adic Galois representation $\rho_{E,p^\infty}:G_F\to G(p^\infty)$ for some $p\mid m$.
    \item \textbf{Entanglement}: The non surjectivity of $\rho_E(G_F)$ on $\underset{\substack{p\text{ prime}\\p\mid m}}{\prod}\rho_{E,p^\infty}(G_F)$.
\end{enumerate}
If $F=\Q$, then $\rho_E$ is never surjective: we have at least that $\rho_{E,2^\infty}$ is not surjective or a Serre entanglement, see \cite[Section 3.13.1]{GroupTheoryEnt}.

Let $F(E[m])$, respectively $F(E[p^\infty])$, be the extension of $F$ generated by the coordinates of the $m$-torsion points of $E$, respectively the $p^k$-torsion points of $E$ for all $k$. The kernel of $\rho_{E,m}$ is $\Gal(\overline{F}/F(E[m]))$. Entanglement is equivalent to the failure of the fields $\{F(E[p^\infty]): p\text{ prime}, p\mid m\}$ to be linearly disjoint. In \cite[Theorem 3.2]{CampagnaStevenhagen} and \cite[Theorem 1.1]{CampagnaPengoEnt}, Campagna, Pengo and Stevenhagen gave, for an elliptic curve $E/F$ with $\End(E)\subset F$, an explicit finite set of primes $S$ such that 

\begin{equation}\label{isoS}
   \Gal(F(E_{\mathrm{tors}})/F)\simeq \Gal(F(E[S^\infty])/F)\times\underset{\substack{p\text{ prime}\\ p\notin S}}{\prod} \Gal(F(E[p^\infty])/F) \tag{$*$}
  \end{equation}
where $F(E[S^\infty])$ is the compositum of the $F(E[p^\infty])$ for $p\in S$. Let $\Delta_F$ be the discriminant of $F$, $\got{f}_E$ be the ideal conductor of $E$ and $\mathrm{N}(\got{f}_E)$ be its norm. If $E/F$ does not have CM, the set $S$ consists of the primes $p$ satisfying at least one of the two following conditions:
\begin{itemize}
    \item $p\mid 2\cdot 3\cdot 5\cdot\Delta_F\cdot \mathrm{N}(\got{f}_E)$,
    \item $\rho_{E,p}$ is not surjective.
\end{itemize}
If $E/F$ has CM by $\mathcal{O}$, an order in $K$, the set $S$ consists of the primes dividing \[[\mathcal{O}_K:\mathcal{O}]\cdot \Delta_F\cdot \mathrm{N}(\got{f}_E).\]

\bigskip

We describe what local conditions and entanglement imply about the division fields $F(E[m])$ for $m\geq2$.

\noindent\textbf{Local conditions.} We know that if $p\geq 5$ and $\rho_{E,p}$ is surjective, then $\rho_{E,p^\infty}$ is surjective, see \cite[Lemma 1]{antwerpIII}. But, if $\rho_{E,p}$ is not surjective, then we cannot deduce the image of $\rho_{E,p^\infty}$ from the image of $\rho_{E,p}$. It can be the full inverse image of $\rho_{E,p}(G_F)$ in $\GL(\Z_p)$, which is equivalent to having \[\Gal(F(E[p^{k+1}])/F(E[p^k]))\simeq(\Z/p\Z)^4\] for all $k\geq1$, but it can be smaller.

\noindent\textbf{Entanglement.} We say that $E/F$ has an $(m,n)$-entanglement if \[F(E[m])\cap F(E[n])\neq F(E[\gcd(m,n)]).\]This corresponds to the linear dependance of the fields $F(E[m])$ and $F(E[n])$ over $F(E[\gcd(m,n)])$.

\noindent Our work focuses on the extreme case: a \textbf{coincidence}. We say that $E/F$ has an \emph{$(m,n)$-coincidence} if $F(E[m])=F(E[n])$. Stevenhagen asked whenever we have a $(2^k,2^{k+1})$-coincidence for an elliptic curve defined over $\Q$ and the answer was given by Rouse and Zureick-Brown in \cite[Remark 1.5]{2adicimage}. Around the same time, in \cite{23Entanglement}, Brau and Jones gave a parametrization of elliptic curves $E/\Q$ such that $\Q(E[2])\subseteq \Q(E[3])$, which is equivalent to having a $(3,6)$-coincidence and corresponds to a $(2,3)$-entanglement. In \cite{coincidences}, Daniels and Lozano-Robledo ask when we have $\Q(\zeta_{p^k})\subseteq\Q(E[m])$ for a prime $p$ and, as a consequence of the Weil pairing, use it to study coincidence. In particular, they showed that the only possible $(m,n)$-coincidence with $m$ and $n$ prime powers are $(2,4)$ and $(2,3)$-coincidences.

In this article, we study the question of coincidence for elliptic curves over an arbitrary number field. As a first observation, we deduce from the isomorphism \eqref{isoS} the following result:

\begin{lemma*}[Lemma {\ref{premiers coincidences dans S}}]Let $m$ and $n$ be two integers and set $m_S$, respectively $n_S$, to be the greatest divisor of $m$, respectively $n$, with only prime divisors in $S$. Suppose that $F(E[m])=F(E[n])$. Then \[F(E[m_S])=F(E[n_S])\]and\[\forall p\notin S,\quad F(E[p^{v_p(m)}])=F(E[p^{v_p(n)}]),\] where $v_p$ denotes the $p$-adic valuation.\end{lemma*}

\noindent We consider the explicit set $S$ given by Campagna et al. If $E/F$ does not have CM, then, for an $(m,n)$-coincidence, we have $v_p(m)=v_p(n)$ for all $p\notin S$. However, the set $S$ is not minimal for this property. In this article, we prove the following result:

\begin{theorem*}[Theorem~\ref{coincidence ramified or bad red}]
Let $E/F$ be an elliptic curve and $m,n\geq1$. Suppose that $F(E[m])=F(E[n])$. Then, for all primes $p$ such that $v_p(m)\neq v_p(n)$, we have $p\mid 2\cdot\Delta_F\cdot\mathrm{N}(\got{f}_E)$.
\end{theorem*}

\begin{corollaire*}[Corollary~\ref{cas F=Q pour coincidence}]
Let $E/\Q$ be an elliptic curve, $\Delta_E$ be the minimal discriminant of $E$ and $m,n\geq1$. Suppose that $\Q(E[m])=\Q(E[n])$. Then, for all primes $p$ such that $v_p(m)\neq v_p(n)$, we have $p\mid 2\Delta_E$.
\end{corollaire*}

If $m\mid n$, then $F(E[m])\subseteq F(E[n])$. For arbitrary positive integers $m$ and $n$, if a coincidence $F(E[m])=F(E[n])$ holds, then it remains true replacing $n$ by $\mathrm{lcm}(m,n)$. Thus, to obtain constraints on coincidences, it suffices to consider $m$ dividing $n$. Furthermore, we can reduce to the question of whether $F(E[m])=F(E[p^km])$ for a prime $p$ and $k\geq1$. Moreover, considering a set $S$ satisfying the isomorphism~\eqref{isoS}, it suffices to consider $m$ with only prime divisors in $S\cup\{p\}$. Then, we consider the following guiding question:

\begin{question}Let $p$ be a prime, $k\geq1$ and $m$ be an integer with only prime divisors in $S\cup\{p\}$. When do we have $F(E[m])=F(E[p^km])$?
\end{question}

\noindent We can reformulate this question, considering $p\nmid m$ and the following situations:
\begin{enumerate}
    \item[-] \emph{Horizontal coincidences} \begin{itemize}
        \item $F(E[m])=F(E[p^km])$ for some $k\geq1$
    \end{itemize}
    \item[-] \emph{Vertical coincidences} \begin{itemize}
        \item $F(E[m])\neq F(E[pm])=\dots=F(E[p^km])$ for some $k\geq2$, or
        \item $F(E[2m])\neq F(E[4m])=\dots=F(E[2^km])$ for some $k\geq3$.
    \end{itemize}
\end{enumerate}
We know by Theorem~\ref{indice coincidence det} that there are no other cases. In addition to refining the set of possible prime divisors as seen above, we give constraints on the possible exponents on the prime divisors. The two following theorems concern horizontal coincidences.

\begin{theorem*}[Corollary~\ref{geatest prime divisor coincidence}]Let $m,n\geq3$, $p$ be a prime such that $p>q$ for all primes $q\mid m$. Suppose that $\Q(\zeta_{p^{v_p(n)}})\not\subseteq F$. If $E/F$ has an $(m,n)$-coincidence, then $v_p(n)=1$.
\end{theorem*}

\begin{theorem*}[Corollary~\ref{coincidence and reduction} and Remark~\ref{coincidence and reduction farther}]
Let $E/\K$ be an elliptic curve with an $(m,n)$-coincidence. Let $p$ be a prime such that $p\mid n$ and $p\nmid m$. Let $\got{p}$ an ideal above $p$ such that the ramification index of $F/\Q$ at $\got{p}$ is prime to $\varphi(p^k)$. Then we are in one of the following cases:
\begin{itemize}
    \item $v_p(n)=1$ and $E/F$ has bad reduction at $\got{p}$,
    \item $v_p(n)=2$, $p=2$ and at $\got{p}$, $E/\K$ has either additive or non split multiplicative reduction,
    \item $v_p(n)=2$, $p=3$, and $E/\K$ has additive and potential good reduction at $\got{p}$,
    \item $v_p(n)=3$ or $4$, $p=2$ and $E/\K$ has additive and potential good reduction at $\got{p}$.
\end{itemize}
\end{theorem*}

We next consider vertical coincidences. The following theorem holds without any hypothesis either on the curve $E/F$ or the field $F$.

\begin{theorem*}[Theorem~\ref{indice coincidence det} and Proposition~\ref{p^4 divise l'indice}]Let $q=p$ and $k\geq1$ if $p$ is odd, or $q=p^2$ and $k\geq2$ if $p$ is even. Let $E/F$ be an elliptic curve. If $F(E[p^k])=F(E[p^{k+1}])$, then $F(E[q])=F(E[p^{k+1}])$ and $p^4\mid [G:\rho_E(G_F)]$.
\end{theorem*}

If $F=\Q$, we know by \cite[Theorem 1.4]{coincidences} that a $(p^k,p^{k+1})$-coincidence is possible only for $p=2$. We prove that, more generally, this is true if $F\cap\Q(\zeta_{p^k})=\Q$, this is Corollary~\ref{vertical coincidence and trivial intersection}. If $E$ has CM by a quadratic field $K$ and $F\subseteq K(j(E))$, then Proposition~\ref{coincidence CM} is more precise: a $(p^k,p^{k+1})$-coincidence is not possible for $p$ odd and $k\geq2$.

Furthermore, we observe that, if $m\mid n$, then the reduction map $\rho_{E,n}(G_F)\to\rho_{E,m}(G_F)$ is an isomorphism. This motivates the study of coincidences in chains \[F(E[m])=F(E[pm])=\dots=F(E[p^km])\]as formulated in the guiding question, and the associated problem of splittings of the surjections $\rho_{E,pm}(G_F)\to\rho_{E,m}(G_F)$ introduced in Subsection~\ref{section split liftable}. We prove the following result:

\begin{theorem*}[Theorem~\ref{T and vertical coincidence}]
Let $p$ be a prime. Let $q=p$ if $p\neq 2,3$, or $q=p^2$ if $p=2,3$. If $\rho_{E,q}(G_F)$ contains a conjugates of $\begin{pmatrix}1&1\\0&1\end{pmatrix}$, then $F(E[q])\neq F(E[pq])$.
\end{theorem*}

The next result is both valid for both horizontal and vertical coincidences. Obviously, an $(m,n)$-coincidence is not possible if $\rho_{E,m}$ and $\rho_{E,n}$ are both surjective, since $\GL(m)$ and $\GL(n)$ are not isomorphic for $m\neq n$. Daniels and Lozano-Robledo compared the abelian part of the division field to show that $E/\Q$ does not have $(m,n)$-coincidence if only $\rho_{E,m}$ is surjective. We use the same idea to show a similar result in case where $\rho_{E,m}$ is large.

\begin{theorem*}[Theorem~\ref{large image m odd}]
Let $m$ be an odd integer and $n\nmid m$ such that $\zeta_n\notin F$.
Suppose that $\rho_{E,m}(G_\K)$ contains $\SL(\Z/m\Z)$ and that $E/F$ has an $(m,n)$-coincidence. Then \[3\mid m,\quad\text{and}\quad \mathrm{D}(\rho_{E,m}(G_F))\neq\SL(\Z/m\Z),\quad \text{and}\quad F(\zeta_n)\subseteq L\] with $L$ a $\Z/3\Z$-extension of $F(\zeta_m)$.
\end{theorem*}

Up to this point, we have only examined necessary condition for having a coincidence. We end the presentation of the result in this article with the statement of sufficient conditions. Over $\Q$, the only possible $(p^k,p^{k+1})$-coincidence for a prime $p$ is a $(2,4)$-coincidence. It is even conjectured that the known $(2,4)$, $(2,3)$, $(2,6)$ and $(3,6)$-coincidences are the only possible coincidences over $\Q$. For any elliptic curve $E/F$, we can construct an $(m,n)$-coincidence with a base change from $F$ to $F(E[\mathrm{lcm}(m,n)])$, but such a base change provides a trivial construction. In Section~\ref{section 3} we prove the existence of a $(4,8)$-coincidence with a base change from $\Q$ to an extension linearly disjoint from $\Q(E[4])$:

\begin{theorem*}[Theorem~\ref{example (4,8) coincidence}]
There are infinitely many isomorphism classes of elliptic curves $E/\Q$ such that there exists a number field $L$ with Galois group $(\Z/2\Z)^r$ over $\Q$ with $1\leq r\leq 4$ satisfying \[L(E[4])=L(E[8])\neq L.\]
\end{theorem*}

\noindent Together with this theorem, in Subsection~\ref{section 3.1} we define the notion of a \emph{minimal base change} and give necessary and sufficient condition to construct an $(m,mn)$-coincidence by a minimal base change of the ground field.

\bigskip

The structure of the paper is as follows. Section~\ref{section 1} contains a useful consequence of the Weil pairing and preliminary constraints for a coincidence. Section~\ref{section 2} deals with horizontal coincidences using the link between the ramification of $F(E[m])/F$ and the reduction type of $E/F$. Section~\ref{section 3} focuses on vertical coincidences. Section~\ref{section 4} concerns the case where $\rho_{E,m}(G_F)$ contains $\SL(\Z/m\Z)$. In Section~\ref{section 5}, we use results of previous sections to give obstructions to a coincidence $F(E[m])=F(E'[n])$ where $E/F$ and $E'/F$ are different elliptic curves.

\section*{Acknowledgement}

I thank my supervisors Samuele Anni and David Kohel for their many useful and interesting suggestions. I also thank Luis Dieulefait for discussion about constructing a non-trivial coincidence over a number field which does not happen over $\Q$.

\section{Notations and preliminaries}\label{section 1}

Through this article we will use the following notations:
\begin{itemize}[label=$-$]
    \item For a finite set $S$, we denote by $\#S$ its cardinality.
    \item For a positive integer $m$ and a prime $p$, $v_p(m)$ denotes the valuation of $m$ at $p$.
    \item For a positive integer $m$, we set \[\GL(m):=\GL(\Z/m\Z)\quad\text{and}\quad \SL(m):=\SL(\Z/m\Z).\]
    \item For finite extensions $K\subseteq L\subseteq M$ such that $M/L$ is Galois, and a prime ideal $\got{p}$ of $L$, we denote by $e_\got{p}(M/K)$ the ramification index of $M/K$ at $\got{p}$. 
    \item $F$ is a number field, $\overline{F}$ denotes an algebraic closure of $F$ and $G_F=\Gal(\overline{F}/F)$ its absolute Galois group.
    \item $(\zeta_m)_{m\geq1}$ is a compatible system of primitive $m$-th root of unity in $\overline{F}$, that is $\zeta_{mn}^n=\zeta_m$ for $m,n\geq1$.
    \item For a prime $p$, $\Q(\mu_{p^\infty})$ denotes the number field generated by the $\zeta_{p^k}$, for $k\geq1$.
    \item $\Delta_F$ denotes the discriminant of the number field $F$.
    \item $\mathcal{O}_F$ denotes the ring of integers of $F$.
\end{itemize}

First of all, we prove Proposition~\ref{premiers coincidences dans S} given in the introduction after some preliminaries. Let $E/F$ be an elliptic curve. If $E/F$ has CM, we suppose that~$F$ contains the CM field. From Serre's open image theorem \cite[Section 4.4, Theorem 3']{Serre} for the non-CM case, and from \cite{CampagnaPengoEnt} for the CM case, there exists a finite set $S$ of rational primes such that \begin{equation}\label{iso entanglement}\Gal(F(E_{\mathrm{tors}})/F)\simeq \Gal(F(E[S^\infty])/F)\times\underset{\substack{p\text{ prime}\\ p\notin S}}{\prod} \Gal(F(E[p^\infty])/F)\end{equation}
where $F(E[S^\infty])$ is the compositum of $F(E[p^\infty])$ for $p\in S$.  The set $S$ is not unique and we are interested in finding the smallest possible set $S$. As underlined in the introduction, Campagna, Pengo and Stevenhagen gave a possible choice of $S$, distinguishing CM and non-CM case. All tensor products are taken over the base field $F$.

\begin{remark}\label{decomp K1K2}Let $I$ be a countable set, $(L_i)_{i\in I}$ be a family of linear disjoint extensions of $F$ and $K/F$ be a subfield of $\underset{i\in I}{\bigotimes} L_i$. Suppose that $K=\underset{i\in I}{\bigotimes}K_i$ with $K_i\subseteq L_i$ for all $i\in I$. For all $i\in I$, we have $K_i\subseteq K\cap L_i$ and \[\underset{i\in I}{\bigotimes} K_i\subseteq \underset{i\in I}{\bigotimes} (K\cap L_i)\subseteq K=\underset{i\in I}{\bigotimes} K_i.\]Hence $K_i=K\cap L_i$. We recall that $\bigotimes K_i$ is a field if and only if the extensions $K_i/F$ are linearly disjoint, which means by definition that the surjection $\bigotimes K_i\to\prod K_i$ on the compositum is an isomorphism.
\end{remark}

We deduce the following proposition, which does not depends on the set $S$. For an integer $m$, let $\gcd(m,S)$ be the greatest divisor of $m$ with only prime divisors in $S$.

\begin{lemma}\label{premiers coincidences dans S}Let $m$ and $n$ be two positive integers and suppose that $E/F$ has an $(m,n)$-coincidence. Then \[F(E[\gcd(m,S)])=F(E[\gcd(n,S)])\]and\[\forall p\notin S,\quad F(E[p^{v_p(m)}])=F(E[p^{v_p(n)}]).\] 
\end{lemma}

\begin{proof}From the isomorphism~\eqref{iso entanglement}, we know that $F(E_\mathrm{tors})=\underset{i\in I}{\bigotimes} F(E[i^\infty])$ for \[I=\{S\}\cup \{p: p \text{ is prime and }p\notin S\}.\]
By linear independance of Galois extensions, we have the following decompositions \[F(E[m])=F(E[\gcd(m,S)])\cdot\underset{p\notin S}{\prod} F(E[p^{v_p(m)}]),\]\[F(E[n])=F(E[\gcd(n,S)])\cdot\underset{p\notin S}{\prod} F(E[p^{v_p(n)}]).\]
By Remark~\ref{decomp K1K2} these decompositions are unique, therefore the proposition holds true.
\end{proof}

We will rely on this consequence of the Weil pairing:

\begin{prop}\label{weil pairing}Let $E/F$ be an elliptic curve and $m$ be a positive integer. The division field $F(E[m])$ contains  $\K(\zeta_m)$. In particular, the image of $\det\circ\rho_{E,m}$ is equal to the image of $\Gal(F(\zeta_m)/F)$ in $(\Z/m\Z)^*$.
\end{prop}

\begin{proof}
The inclusion $F(\zeta_m)\subseteq F(E[m])$ follows from the Galois invariance of the Weil pairing, see \cite[Corollary 8.1.1]{AEC}. Let $e_m$ be the Weil pairing on $E[m]$ and $(P,Q)$ be a basis of the $m$-torsion such that $e_m(P,Q)=\zeta_m$.  Now, for $\sigma\in \Gal(F(E[m])/F)$, the Galois invariance of $e_m$ gives 
\[e_m(\rho_{E,m}(\sigma)(P,Q))=e_m(\sigma(P),\sigma(Q))=\sigma(e_m(P,Q))=\sigma(\zeta_m),\]
and the Weil pairing being bilinear and alternating, we have
\[e_m(\rho_{E,m}(\sigma)(P,Q))=e_m(P,Q)^{\det\circ\rho_{E,m}(\sigma)}=\zeta_m^{\det\circ\rho_{E,m}(\sigma)}.\]
The result follows.
\end{proof}

Proposition~\ref{weil pairing} implies that, if $F(E[n])=F(E[m])$ for some $n\geq m$, then $F(\zeta_n)\subseteq F(E[m])$. A recurring strategy will be to give restrictions on having this inclusion.

\begin{remark}
Let $n$ and $m$ be two integers such that $m<n$. Then there exists a prime $p$ such that $v_p(n)-v_p(m)=k\geq1$. On the one hand, we have \[F(E[m])\subseteq F(E[p^km])\subseteq F(E[\mathrm{lcm}(m,n)])=F(E[m])F(E[n]).\]On the other hand, we have $F(\zeta_{p^km})\subseteq F(E[p^km])$. Therefore, each time we have $F(\zeta_{p^km})\not\subseteq F(E[m])$ for some $k\geq1$, it follows that $F(E[m])\neq F(E[n])$ for all $n$ such that $v_p(n)-v_p(m)\geq k$.
\end{remark}

As a first attempt, we investigate the possibility of an $(m,n)$-coincidence, simply by using the inclusions of fields $F(\zeta_m)\subseteq F(E[m])$ and of groups $\rho_{E,m}(G_F)\leq \GL(m)$, and the resulting divisibility of degrees and orders.  The next proposition tells us that, if $F(E[n])=F(E[m])$, subject to an additional condition on $F$, then the primes greater than every prime dividing $m$ can divide $n$ to at most power $1$, unless $m$ is a power of $2$, in which case $3$ can divide $n$ with possibly a greater power than $1$.

\begin{prop}\label{greatest prime divisor}Let $m\geq2$, $p$ be a prime such that $p>q$ for all primes $q\mid m$ and $r$ be the largest integer such that  $\Q(\zeta_{p^{r}})\subseteq \K\cap\Q(\mu_{p^\infty})$.
Let $E/\K$ be an elliptic curve such that $\K(\zeta_{p^k})\subseteq \K(E[m])$ with $k>r$. Then, $k=1$ (and $r=0$), unless $(m,p)=(2^j,3)$ for some $j\geq1$, in which case either $r=0$ and $k\leq 2$, or $r=k-1$.
\end{prop}

\begin{proof}Suppose that $r>0$, or $r=0$ and $k\geq2$. We will prove that $(m,p)=(2^j,3)$, and $r=k-1$ or $r=0$ and $k=2$.
By assumption, $p^{k-r}\mid[\K(\zeta_{p^k}):\K]$ if $r>0$ or $p^{k-1}\mid [\K(\zeta_{p^k}):\K]$ if $r=0$. In any case, $p$ divides $[F(\zeta_{p^k}):F]$. Since $\K(\zeta_{p^k})\subseteq \K(E[m])$, we have
\begin{equation*}[\K(\zeta_{p^k}):\K]\bigm| [\K(E[m]):\K]\bigm| \# \GL(m),\end{equation*}and \[\#\GL(m)=\underset{\substack{q^j\mid m\\j=v_q(m)}}{\prod}\# \GL(q^j)=\underset{\substack{q^j\mid m\\ j=v_q(m)}}{\prod}q^{4(j-1)+1}(q-1)^2(q+1).\]
Therefore, since $q<p$ for all $q\mid m$, we obtain $p=q+1$ for some $q$ dividing $m$ and so $m=2^j$ for some $j\geq1$ and $p=3$. In this case, $k-r=1$ if $r>0$ and $k-1=1$ if $r=0$.
\end{proof}

\begin{corollaire}\label{geatest prime divisor coincidence}Under the hypotheses of Proposition~\ref{greatest prime divisor}, let $n$ be an integer such that $v_p(n)=k$ and suppose that $E/F$ has an $(m,n)$-coincidence. Then, $k=1$, unless $(m,p)=(2^j,3)$ for some $j\geq1$, in which case $r=0$ and $k\leq 2$, or $r=k-1$.
\end{corollaire}

In Corollary~\ref{coincidence and reduction}, in the next section, we extend Proposition~\ref{greatest prime divisor} by replacing "$q<p$ for all $q\mid m$" by "$p\nmid m$", at the expense of adding conditions on the ramification at $p$ in $F$ or on the reduction type of $E$ at $p$. 


\section{Horizontal coincidence : ramification behaviour}\label{section 2}

We talk about $\emph{horizontal coincidence}$ if we have an $(m,n)$-coincidence and the sets of prime divisors of $m$ and $n$ are not the same. In this section, we study the obstructions to horizontal coincidences given by the type of reduction of the elliptic curve and the resulting ramification. 

\subsection{Ramification and reduction type}

Let $\got{p}$ be a prime of $\mathcal{O}_F$, whose residue characteristic is $p$.
We recall the criterion of Néron-Ogg-Shafarevitch:
\begin{prop}[{\cite[VII, Theorem 7.1]{AEC}}]\label{Criterion NOS}
Let $E/\K$ be an elliptic curve. If $E/\K$ has good reduction at $\got{p}$, then $\K(E[m])/\K$ is unramified at $\got{p}$ for all $m$ such that $p\nmid m$.
\end{prop}

Moreover, the theory of Tate curves gives constraints on the ramification when the reduction is multiplicative:

\begin{prop}\label{rami red mult}
Let $E/F$ be an elliptic curve and $m\geq2$ such that $p\nmid m$. If $E/\K$ has split multiplicative reduction at $\got{p}$ or if $E/\K$ has multiplicative reduction $\got{p}$ and $p$ is odd, then $\K(E[m])/\K$ is tamely ramified at $\got{p}$. If $E/F$ has non split multiplication at $\got{p}$ and $p$ is even, then $v_p(e_\got{p}(F(E[m])/F))\leq1$.
\end{prop}

\begin{proof}First, we suppose that $E/F$ has split multiplicative reduction at $\got{p}$.
Let $F_\got{p}$ be the completion of $F$ at $\got{p}$. We have, from \cite[V.Theorem 5.3]{ATAEC}, that $E$ is isomorphic over $F_\got{p}$ to the Tate curve $E_q$ for some $q\in \K_\got{p}^*$ (for the definition of $E_q$, see \cite[V.Theorem 3.1]{ATAEC}). We consider the $\got{p}$-adic uniformization: \[\overline{\K_\got{p}}^*/q^\Z\overset{\sim}{\longrightarrow}E_q(\overline{\K_\got{p}}).\]
Restricting to the group of $m$-torsion on each side, we obtain an isomorphism \[\phi:\left(\zeta_m^\Z Q^\Z\right)/q^\Z\overset{\sim}{\longrightarrow}E_q[m],\]where $Q=q^\frac{1}{m}$ is a $m$-th root of $q$. 
The action of $\Gal(\overline{\K_\got{p}}/F_\got{p})$ on $E_q[m]$ is compatible with its action on $\left(\zeta_m^\Z Q^\Z\right)/q^\Z$ (see \cite[V, Theorem 5.3]{AEC}). Let $I_\got{p}$ be the inertia subgroup of $\Gal(\overline{\K_\got{p}}/F_\got{p})$ and let $\sigma\in I_\got{p}$. Since $p\nmid m$, the extension $\K_\got{p}(\zeta_m)/F_\got{p}$ is unramified, and so $\sigma(\zeta_m)=\zeta_m$. Since $Q$ is a root of $X^m-q$, so is $\sigma(Q)$. Therefore, there exists $a\in\Z/m\Z$ such that $\sigma(Q)=\zeta_m^aQ$. We set $P_1=\phi(\zeta_m)$ and $P_2=\phi(Q)$. Then \[\sigma(P_1)=\sigma(\phi(\zeta_m))=\phi(\sigma(\zeta_m))=\phi(\zeta_m)=P_1\]and, \[\sigma(P_2)=\sigma(\phi(Q))=\phi(\sigma(Q))=\phi(\zeta_m^a Q)=a\phi(\zeta_m)+\phi(Q)=aP_1+P_2.\]
Hence, for all $\sigma\in I_\got{p}$, there exists $a\in\Z/m\Z$ such that \[\rho_{E,p}(\sigma)=\begin{pmatrix}1&a\\0&1\end{pmatrix}.\]
It follows that the image of the wild inertia by $\rho_{E,p}$ is included in a group of order $m$. However, as observed by Serre in \cite[Section 1.1]{Serre}, $I_\got{p}$ is a pro-$p$-group, and so its image by $\rho_{E,p}$ is a $p$-group. So it is trivial. Hence, $\K(E[m])/\K$ is tamely ramified at $\got{p}$.

Now, suppose that $E/\K$ has non split multiplicative reduction at $\got{p}$, and let $L/F$ be the quadratic extension where the reduction is split. Then $L(E[m])/L$ is tamely ramified at $p$. Moreover, \[e_\got{p}(L(E[m])/F)=e_\got{p}(L(E[m])/L)e_\got{p}(L/F)\] and so 
\[e_\got{p}(F(E[m])/F)\mid e_\got{p}(L(E[m])/F)\mid 2e_\got{p}(L(E[m])/L),\]
which completes the proof.
\end{proof}

Finally, we also have constraints on the ramification in case of additive reduction: 

\begin{prop}\label{rami additive reduction}
Let $E/\K$ be an elliptic curve and $m\geq2$ such that $p\nmid m$. If $p>3$ and $E/\K$ has additive reduction at $\got{p}$, or if $p=3$ and $E/F$ does not have potential good reduction at $\got{p}$, then $F(E[m])/F$ is tamely ramified at $\got{p}$.
\end{prop}

\begin{proof}
If $E/F$ has potential good reduction, the proposition follows from \cite[Section 2, Corollary 2]{goodredabvar}. If $E/F$ does not have potential good reduction, then the results follows from Proposition~\ref{rami red mult} and \cite[Appendix C, Theorem 14.1]{AEC}, since a quadratic extension cannot be widely ramified outside $2$.
\end{proof}

Finally, let us recall the following result:

\begin{prop}[{\cite[Section 4.2]{samuelelocalglobal}}]\label{stable red deg 24}
Let $E/\K$ be an elliptic curve and $m\geq2$ such that $p\nmid m$. Suppose that $E/F$ has additive reduction at $\got{p}$. There exists an extension $L/F$ of degree dividing $24$ such that $E/L$ has stable reduction at $\got{p}$.
\end{prop}

\subsection{Ramification and entanglement}

Let $p$ be a prime and $\got{p}$ be a prime ideal of $F$ above $p$. Set $e= e_\got{p}(F/\Q)$ the ramification index of $\got{p}$ in $F/\Q$. We know that, if $F(E[n])\subseteq F(E[m])$ then $F(\zeta_{p^k})\subseteq F(E[m])$ for all $p^k\mid n$. In particular, $e_\got{p}(F(\zeta_{p^k})/F)$ divides $ e_\got{p}(F(E[m])/F)$ for all $p^k\mid n$. Lemma~\ref{rami ext cyclo} gives information about $e_\got{p}(F(\zeta_{p^k})/F)$.

The map $\varphi:\Z\to\Z$ denotes the Euler totient function. 

\begin{lemma}\label{rami ext cyclo}
We have $v_p(e)\geq k-1-v_p(e_\got{p}(F(\zeta_{p^k})/F))$. Moreover, if $e_\got{p}(F(\zeta_{p^k})/F)=1$, then $\varphi(p^k)\mid e$.
\end{lemma}

\begin{proof}
The extension $F(\zeta_{p^k})/F$ is Galois and so the ramification index above $\got{p}$ only depends on $\got{p}$. We have \begin{align*}e_\got{p}(F(\zeta_{p^k})/F)e&=e_\got{p}(F(\zeta_{p^k})/\Q)\\&=e_\got{p}(F(\zeta_{p^k})/\Q(\zeta_{p^{k}}))e_\got{p}(\Q(\zeta_{p^{k}})/\Q)\\&=e_\got{p}(F(\zeta_{p^k})/\Q(\zeta_{p^k}))\varphi(p^k)\cdot\end{align*}
Since $v_p(\varphi(p^k))=k-1$ we obtain the first statement. The second follows from the previous equality.
\end{proof}

The previous section gives information about $e_\got{p}(F(E[m])/F)$, summarized in Theorem~\ref{table of ramification division}. These results give constraints on having an $(m,n)$-coincidence when $m$ and $n$ do not have the same prime divisors.


\begin{theorem}\label{table of ramification division}Let $m\geq2$ such that $p\nmid m$. Let $E/F$ be an elliptic curve. The valuation at $\got{p}$ of the ramification index $e_\got{p}(F(E[m])/F)$ appears in the table below together with sufficient conditions on the reduction of $E/F$ at $\got{p}$.

\vspace*{0.5cm}
\centering{\scalebox{0.9}{\begin{tabular}{|c|c|c|}
  \hline
  Sufficient condition on $E/F$ & $t=e_\got{p}(F(E[m])/F)$\\
  \hline
  good reduction at $\got{p}$ & $t=1$ \\
  \hline
  multiplicative red. at $\got{p}$ with $p$ odd & \\
  split multiplicative red. at $\got{p}$ with $p=2$ & $v_p(t)=0$\\
  additive red. at $\got{p}$ with $p>3$ & \\
  additive, not potentially good red. at $\got{p}$ with $p=3$ & \\
  \hline
  (non split) multiplicative red. at $\got{p}$ with $p=2$ &$v_p(t)\leq1$\\
  additive, potentially good red. at $\got{p}$ with $p=3$ & \\
  \hline
  additive red. at $\got{p}$ with $p=2$ &$v_p(t)\leq 3$\\
  \hline
\end{tabular}}}

\justifying

\vspace*{0.5cm}
\end{theorem}

\begin{proof}
Here is the table, with an additional column with the propositions required for the proof.

\vspace*{0.5cm}

\centering{\scalebox{0.9}{\begin{tabular}{|c|c|c|}
  \hline
  Sufficient condition on $E/F$ & $t=e_\got{p}(F(E[m])/F)$&Proof\\
  \hline
  good red. at $\got{p}$ & $t=1$ & Proposition~\ref{Criterion NOS}\\
  \hline
  mult. red. at $\got{p}$ with $p$ odd & & Proposition~\ref{rami red mult} \\
  split mult. red. at $\got{p}$ with $p=2$ & $v_p(t)=0$ & Proposition~\ref{rami red mult}\\
  add. red. at $\got{p}$, $p>3$ &  & Proposition~\ref{rami additive reduction}\\
  add., no pot. good red. at $\got{p}$ with $p=3$ & & Proposition~\ref{rami additive reduction}\\
  \hline
  (non split) mult. red. at $\got{p}$ with $p=2$ &$v_p(t)\leq1$& Proposition~\ref{rami red mult}\\
  add., pot. good red. at $\got{p}$ with $p=3$ & & Proposition~\ref{stable red deg 24}\\
  \hline
  add. red. at $\got{p}$ with $p=2$ &$v_p(t)\leq 3$& Proposition~\ref{stable red deg 24}\\
  \hline
\end{tabular}}}
\vspace*{0.5cm}

\justifying
\end{proof}

\begin{remark}\label{table of ramification cyclo}By Lemma~\ref{rami ext cyclo}, we obtain the table below, in which we present the necessary condition on the ramification of $F/\Q$ to obtain the ramification index as in previous theorem.

\vspace*{0.5cm}
\centering
\scalebox{0.82}{\begin{tabular}{|c|c|c|}
  \hline
  $s=e_\got{p}(F(\zeta_{p^k})/F)$ & Necessary condition on $F/\Q$ \\
  \hline
   $s=1$& $\varphi(p^k)\mid e$\\
  \hline
   $v_p(s)=0$&$v_p(e)\geq k-1$ \\
  \hline
  $v_p(s)\leq1$&$v_p(e)\geq k-2$\\
  \hline
  $v_p(s)\leq 3$&$v_p(e)\geq k-4$\\
  \hline
\end{tabular}}

\justifying

\vspace*{0.5cm}
\end{remark}

With the notation of Theorem~\ref{table of ramification division} and Remark~\ref{table of ramification cyclo}, if $F(\zeta_{p^k})\subseteq F(E[m])$, then we must have $s\mid t$. Therefore, the tables give restrictions on having $F(\zeta_{p^k})\subseteq F(E[m])$. For example, if we have this inclusion and $E/F$ has good reduction at $\got{p}$, then $\varphi(p^k)$ must divide the ramification index of $F/\Q$ at $\got{p}$. In the following corollary, we consider the case of $F/\Q$ unramified above $p$.

\begin{corollaire}\label{cyclo inclusion and reduction}
Let $E/\K$ be an elliptic curve, $m\geq2$, $k\geq1$ and suppose that $p\nmid m\Delta_F$. If $\K(\zeta_{p^k})\subseteq \K(E[m])$, then we are in one of the following cases:
\begin{itemize}
    \item $k=1$ and $E/F$ has bad reduction at every ideal above $p$,
    \item $k=2$, $p=2$ and at each prime above $p$, $E/\K$ has either additive or non split multiplicative reduction,
    \item $k=2$, $p=3$, and $E/\K$ has additive and potential good reduction at every ideal above $p$,
    \item $k=3$ or $4$, $p=2$ and $E/\K$ has additive and potential good reduction at every ideal above $p$.
\end{itemize}
\end{corollaire}

\begin{proof}Let $\got{p}$ be a prime ideal above $p$. Since $p\nmid\Delta_F$, the extension $F(\zeta_{p^k})/F$ is ramified at $\got{p}$, from Lemma~\ref{rami ext cyclo}, so is $F(E[m])/F$ by assumptions. Looking at the tables in Theorem~\ref{table of ramification division} and Remark~\ref{table of ramification cyclo}, with $e=1$, we see that the possibilities are: $k=1$, corresponding to the second line of each tables; $k=2$, corresponding to the third line and in this case $p=2,3$; or $k=3,4$, corresponding to the fourth line, where $p=2$.
\end{proof}

The corollary below tells us that if $E/F$ has a $(m,n)$-coincidence (with some conditions on $F$), then the primes greater than $5$ not dividing $m$ (respectively~$n$) divide $n$ (respectively $m$) to at most power $1$, and $E/F$ must have bad reduction at these primes. Moreover, if $3\nmid m$, then $3$ divides $n$ to at most power $2$, and if $m$ is odd, then $2$ divides $n$ to at most power $4$, and the greater the power, the more restrictive is the reduction type.

\begin{corollaire}\label{coincidence and reduction}
Let $E/\K$ be an elliptic curve with an $(m,n)$-coincidence. Suppose that $p\mid n$ and $p\nmid m\Delta_F$. Then we are in one of the following cases:
\begin{itemize}
    \item $v_p(n)=1$ and $E/F$ has bad reduction at every ideal above $p$,
    \item $v_p(n)=2$, $p=2$ and at each prime above $p$, $E/\K$ has either additive or non split multiplicative reduction,
    \item $v_p(n)=2$, $p=3$, and $E/\K$ has additive and potential good reduction at every ideal above $p$,
    \item $v_p(n)=3$ or $4$, $p=2$ and $E/\K$ has additive and potential good reduction at every ideal above $p$.
\end{itemize}
\end{corollaire}

\begin{remark}\label{coincidence and reduction farther}If $e_\got{p}(F/\Q)$ is prime to $\varphi(p^k)$ (hypothesis that is satisfied for example if $F/\Q$ is unramified at $\got{p}$), then Corollaries~\ref{cyclo inclusion and reduction} and~\ref{coincidence and reduction} are true replacing "at every ideal above $p$" by "at $\got{p}$".
\end{remark}

Again using ramification, we also deduce a result on vertical coincidences, which is the topic of the next section:

\begin{theorem}\label{ramification et coincidence verticale}Let $E/F$ be an elliptic curve, $p$ be a prime and $k\geq2$. If $p^{k-1}\nmid e_\got{p}(F/\Q)$ and $E/F$ has good supersingular reduction at $\got{p}$, then $F(E[p])\neq F(E[p^k])$.
\end{theorem}

\begin{proof}
From \cite[Proposition 12]{Serre}, the extension $F(E[p])/F$ is tamely ramified at $p$. Since $p^{k-1}$ does not divide $e_\got{p}(F/\Q)$, then Remark~\ref{table of ramification cyclo} gives that $F(E(\zeta_{p^2})/F$ is wildly ramified at $p$. The result follows as a consequence of the Weil pairing~\ref{weil pairing}. 
\end{proof}

\section{Coincidences in towers}\label{section 3}

In this section, we deal with \emph{coincidence in towers}, or \emph{vertical coincidences}, that is to say $(p^k,p^{k+1})$-coincidences for a prime $p$ and a positive integer~$k$. More generally, the section also contains results about $(m,n)$-coincidence where $m\mid n$.

\subsection{Construction of vertical coincidences}\label{section 3.1}

Over $\Q$, we know that infinitely many elliptic curves have a $(2,4)$-coincidence and this is the only vertical coincidence which occurs, see \cite[Theorem 1.4]{coincidences}. Over a number field, there are additional possibilities. Obviously, to obtain an $(m,mn)$-coincidence for an elliptic curve $E/F$, it suffices to do a base change of the ground field to~$F(E[mn])$. However, such a base change is a trivial construction and so not very relevant. We will say that the base change from $F$ to $L$ is \emph{minimal} for an $(m,mn)$-coincidence if $L(E[m])=F(E[mn])$ and $F(E[m])\cap L=F$. Here is an example of a $(4,8)$-coincidence obtained by a minimal base change:

\begin{theorem}\label{example (4,8) coincidence}
There are infinitely many isomorphism classes of elliptic curves $E/\Q$ such that there exists a number field $L$ with Galois group $(\Z/2\Z)^r$ over~$\Q$ with $1\leq r\leq 4$ satisfying \[L(E[4])=L(E[8])\neq L.\]
\end{theorem}

\begin{proof}
We apply Proposition~\ref{constructing coincidence from split} for $F=\Q$, $m=4$, $r=2$ and $E/\Q$ such that $\Gal(\Q(E[8])/\Q)\simeq (\Z/2\Z)^t$ for some $t$. Hence, there exists $L/\Q$ of degree dividing $\#\GL(8)/\#\GL(4)=2^4$ (by \eqref{valeurs des suites d'index} in Subsection~\ref{index of images}) such that $L(E[8])=L(E[4])\neq L$. By \cite[Theorem 1.1]{abdivfield}, the Galois group $\Gal(\Q(E[8])/\Q)$ is isomorphic to $(\Z/2\Z)^t$ for some $t$ for infinitely many isomorphism classes of elliptic curves $E/\Q$ and $\Q(E[4])/\Q$ is non trivial since it contains $\zeta_4$, which completes the proof. 
\end{proof}

\begin{remark}
The proof of Theorem~\ref{example (4,8) coincidence} considers only elliptic curves $E/\Q$ such that $\Q(E[8])/\Q$ is abelian. In this case $\Gal(\Q(E[8])/\Q)\simeq(\Z/2\Z)^t$ with $t\in\{4,5,6\}$ from~\cite[Theorem 1.1]{abdivfield}. Let $r$ be as in Theorem~\ref{example (4,8) coincidence}. We have $2^r=\#\Gal(\Q(E[8])/\Q(E[4]))$ by construction of $F$ and $1\leq r\leq 4$.  If $t=4$, then $1\leq r\leq 3$ and if $t=6$ then $2\leq r\leq 4$.\end{remark}

\begin{remark}
We cannot use the abelian case to construct $(p,p^2)$-coincidence with $p$ odd, because there is no abelian $p^2$-division field for $p$ odd and $E$ defined over $\Q$. 
\end{remark}

\begin{prop}\label{constructing coincidence from split}Let $m,n$ be positive integers.
Let $E/F$ be an elliptic curve such that the following exact sequence is split:\[1\to\Gal(F(E[mn])/F(E[m]))\to\Gal(F(E[mn])/F)\to\Gal(F(E[m])/F)\to1\]with $F(E[m])/F$ non trivial.
 Then there exists an extension $L/F$ of degree dividing $\#\GL(mn)/\#\GL(m)$ such that \[L(E[m])=L(E[mn])\neq L.\]
\end{prop}

To prove the proposition, we will use the following elementary remark:

\begin{remark}\label{divisibilité des indices}Let $G$ be a group and $H$ be a subgroup of $G$ of finite index. Let $\phi:G\to G'$ be a surjective morphism (of groups) and set $H'=\phi(H)$. Then $\phi$ induces a surjective morphism of $G$-sets $G/H\to G'/H'$, from which \[[G:H]=[G':H'][\ker(\phi)H:\ker(\phi)].\]In particular $[G':H']$ divides $[G:H]$.
\end{remark}

\begin{proof}[Proof of Proposition~\ref{constructing coincidence from split}]
Since the sequence is split, there exists a morphism \[\iota:\Gal(F(E[m])/F)\to \Gal(F(E[mn])/F)\] such that the composition with the restriction map \[\Gal(F(E[mn])/F)\to \Gal(F(E[m])/F)\] is the identity. Let $L$ be the fixed field of $\mathrm{Im}(\iota)$. Then $\Gal(F(E[mn])/F)$ is the semi-direct product of $\Gal(F(E[mn])/F(E[m]))$ by $\Gal(F(E[mn])/L)$ and so $L(E[m])=F(E[mn])=L(E[mn])$ and $L\cap F(E[m])=F$. Since $F(E[m])/F$ is nontrivial, then $L(E[m])\neq L$. Moreover, the extension $L/F$ has degree $[F(E[mn]):F(E[m])]$, which divides $\#\GL(mn)/\#\GL(m)$ by point (1) of Remark~\ref{divisibilité des indices}, taking for $\phi$ the natural map $\GL(mn)\to\GL(m)$ and $H=\rho_{E,mn}(G_F)$.
\end{proof}

\begin{corollaire}
For $E/F$ an elliptic curve, the following are equivalent: \begin{enumerate}
    \item The following sequence is split \[\scalebox{0.95}{$1\to\Gal(F(E[mn])/F(E[m]))\to\Gal(F(E[mn])/F)\to\Gal(F(E[m])/F)\to1.$}\]
    \item There exists an injective morphism \[\iota:\Gal(F(E[m])/F)\to\Gal(F(E[mn])/F)\] such that \[\Gal(F(E[mn])/F)=\iota(\Gal(F(E[m])/F))\ltimes \Gal(F(E[mn])/F(E[m])).\]
    \item There exists a minimal base change $L/F$ such that $E/L$ has an $(m,mn)$-coincidence.
\end{enumerate}
In this case, $\iota(\Gal(F(E[m])/F)=\Gal(F(E[mn])/L)$. 
\end{corollaire}

\begin{proof}
The equivalence between point $(1)$ and $(2)$ is immediate from the definitions of split exact sequence and semi-direct product. We know that $(1)\implies (3)$ by Proposition~\ref{constructing coincidence from split}. It remains to show that $(3)\implies (1)$. Suppose that the conditions of $(3)$ are satisfied. Since $F(E[mn])=L(E[mn])$, we have the following commutative diagram, where the horizontal arrows are restriction morphisms and the vertical arrows are inclusion morphisms:
\[\begin{tikzcd}\Gal(F(E[mn])/L)\ar[r, "\psi"]\ar[d, hook]&\Gal(L(E[m])/L)\ar[d, hook, "\phi"]\\\Gal(F(E[mn])/F)\ar[r]&\Gal(F(E[m])/F).\end{tikzcd}\] By assumption, $\psi$ is a isomorphism, together with $\phi$ by linear independance of $F(E[m])$ and $L$ over $F$. It follows that the exact sequence of $(1)$ splits by the morphism
\[\begin{array}{ccccc}
    \iota&:&\Gal(F(E[m])/F)&\longrightarrow& \Gal(F(E[mn])/F)\\ &&\sigma&\longmapsto& (\phi\circ \psi)^{-1}(\sigma).\end{array}\]
\end{proof}

As a consequence of the corollary, the elliptic curves satisfying Theorem~\ref{example (4,8) coincidence} are exactly those such that we have a split exact sequence \[1\to(\Z/2\Z)^r\to\Gal(\Q(E[8])/\Q)\to\Gal(\Q(E[4])/\Q)\to1.\] 
In particular, this is true for $E/\Q$ such that $\Q(E[8])/\Q$ is a $(\Z/2\Z)^t$-extension and a classification for such elliptic curves is given in \cite[Table 4]{abdivfield}. But there are many other possibilities. More generally, we have \[\Gal(F(E[p^{k+1}])/F(E[p^k]))\simeq(\Z/p\Z)^r\] for some $r\leq 4$. Indeed, the Galois group $\Gal(F(E[p^{k+1}])/F(E[p^k]))$ is isomorphic, for $n=2$, to a subgroup of
\begin{equation}\label{noyau de la projection} \ker\left(\mathrm{GL}_n(p^{k+1})\to\mathrm{GL}_n(p^k)\right)=I_n+p^k\mathrm{M}_n(\Z/p^{k+1}\Z)\simeq (\Z/p\Z)^{n^2}.\end{equation} Hence, to construct a $(p^k,p^{k+1})$-coincidence by minimal base change, we have and it suffices to find elliptic curves $E/F$ such that the following exact sequence is split: \[1\to(\Z/p\Z)^r\to\Gal(F(E[p^{k+1}])/F)\to\Gal(F(E[p^k])/F)\to1.\]

\subsection{Trivial intersection with the cyclotomic field}\label{section 3.2}

In this section, we show that, if $F\cap\Q(\zeta_{p^k})$ is trivial, then a $(p^k,p^{k+1})$-coincidence with $p$ prime is possible only for $p=2$.

\begin{lemma}\label{2.5}Let $E/\K$ be an elliptic curve, and $L/F$ be a cyclic extension such that $L\subseteq  \K(E[m])$. Let $\sigma\in G_\K$ such that its restriction to $L$ generates $\Gal(L/\K)$. Then the order of $\rho_{E,m}(\sigma)$ is divisible by $[L:\K]$.\end{lemma}


\begin{proof}Let $\barr{\rho_{E,m}}$ be the reduction of $\rho_{E,m}$ modulo $\Gal(\overline{F}/F(E[m]))$. Then \[[L:\K]=\ord(\restriction{\sigma}{L})\mid\ord(\restriction{\sigma}{F(E[m])})=\ord(\barr{\rho_{E,m}}(\restriction{\sigma}{F(E[m])}))=\ord(\rho_{E,m}(\sigma)).\]The first equality is by assumption, and the second is because of the injectivity of $\overline{\rho_{E,m}}$.\end{proof}


\begin{theorem}\label{equality for p odd}
Let $E/\K$ be an elliptic curve, $p$ be a prime and $k$ be a positive integer such that $\K\cap\Q(\zeta_{p^k})=\Q$. If $\K(\zeta_{p^{k+1}})\subseteq \K(E[p^k])$, then $p=2$.
\end{theorem}

\begin{proof}Suppose that $p$ is odd and $\K(\zeta_{p^{k+1}})\subseteq \K(E[p^k])$. Let $\sigma\in G_\K$ such that its restriction to $\K(\zeta_{p^{k+1}})$ generates $\mathrm{Gal}(\K(\zeta_{p^{k+1}})/\K)$. Then its restriction to $\K(\zeta_{p^k})$ generates $\Gal(\K(\zeta_{p^k})/\K)$. So $\det\rho_{E,p^k}(\sigma)$ generates $(\Z/p^k\Z)^*$. Moreover, Lemma~\ref{2.5} says that $\varphi(p^{k+1})$ divides the order of $\rho_{E,p^k}(\sigma)$ and so its determinant is a square mod $p$, by \cite[Lemma 3.5]{coincidences}. But, for $p$ odd, a square mod $p$ cannot generate $(\Z/p^k\Z)^*$. Hence, $p$ is even.
\end{proof}

\begin{remark}\label{vertically2ent}
The CM elliptic curve $y^2=x^3-11x-14$ satisfies $\Q(\zeta_{2^{k+1}})\subseteq\Q(E[2^k])$ for all $k\geq1$. See \cite[Theorem 1.5]{coincidences}, and \cite[Theorem 1.1]{vertically2ent} for more examples of such curves.
\end{remark}

\begin{corollaire}\label{vertical coincidence and trivial intersection}
Let $E/\K$ be an elliptic curve with $\K\cap\Q(\zeta_{p^k})=\Q$. If $\K(E[p^k])=\K(E[p^{k+1}])$, then $p=2$.
\end{corollaire}

\begin{proof}
It is immediate from Theorem~\ref{equality for p odd}, since $\K(\zeta_{p^{k+1}})\subseteq \K(E[p^{k+1}])$.
\end{proof}

\begin{remark}
In \cite[Theorem 1.4]{coincidences}, Daniels and Lozano-Robledo have shown that, for $F=\Q$, only $k=1$ occurs.
\end{remark}

\begin{remark} If $E/F$ has an $(m,mn)$-coincidence, then we must have the inclusion $F(\zeta_{mn})\subseteq F(E[m])$. But this does not implies in general that $F(\zeta_{mn})=F(\zeta_{m})$, as we will see in Remark~\ref{suite ell_k}. Even more, unless $m$ is odd and $n=2$, this last never happens if $F=\Q$, and yet some coincidences occurs, like $(2,4)$ and $(2,6)$-coincidence, see \cite[Examples 1.2 and 1.3]{coincidences}. As in Remark~\ref{suite ell_k}, it is due to the non-surjectivity of $\SL(mn)\cap \rho_{E,mn}(G_F)\to \SL(m)\cap \rho_{E,m}(G_F)$. 
\end{remark}

\begin{remark}The condition $F(\zeta_{p^{k+1}})\subseteq F(E[p^k])$ is not sufficient to have the coincidence. For example, the elliptic curve of Remark~\ref{vertically2ent} does not satisfy $\Q(E[2^k])=\Q(E[2^{k+1}])$ for any $k$. Indeed, this elliptic curve has CM and \cite[Proposition 3.9]{coincidences} implies that no CM elliptic curve defined over $\Q$ has a $(2^k,2^{k+1})$-coincidence. 
\end{remark}

We are now able to prove the theorem stated in the introduction:

\begin{theorem}\label{coincidence ramified or bad red}
Let $m,n\geq1$ and $E/F$ be an elliptic curve with conductor ideal $\got{f}_E$. Let $\mathrm{N}(\got{f}_E)$ be the norm of $\got{f}_E$. Suppose that $F(E[m])=F(E[n])$. Then, for all primes $p$ such that $v_p(m)\neq v_p(n)$, we have \[p\mid2\cdot\Delta_F \cdot\mathrm{N}(\got{f}_E).\] 
\end{theorem}

\begin{proof}
First, suppose that $p$ divides $n$ or $m$ but not both. Then, by Corollary~\ref{coincidence and reduction}, if $p\nmid\Delta_F$, then $E/F$ has bad reduction above $p$.
Now, suppose that $p$ divides both $n$ and $m$ such that $v_p(m)=k$ and $v_p(m)<v_p(n)$. Since $F(E[m])=F(E[\mathrm{lcm}(m,n)])$, then $F(E[m])=F(E[pm])$. Setting $a=\frac{m}{p^k}$ and $L=F(E[a])$, we obtain $F(E[p^ka])=F(E[p^{k+1}a])$ and $L(E[p^k])=L(E[p^{k+1}])$. Then Corollary~\ref{vertical coincidence and trivial intersection} implies that $L\cap \Q(\zeta_{p^{k+1}})\neq\Q$ or $p=2$. In particular, $p$ is ramified in $L/\Q$ or $p=2$. But $L=F(E[a])$, so $p$ is ramified in $L/\Q$ if and only if $p$ is ramified in $F/\Q$ or in $F(E[a])/F$. Therefore $p\mid\Delta_F$ or $E$ has bad reduction above $p$.
\end{proof}

\begin{corollaire}\label{cas F=Q pour coincidence}
Let $m,n\geq1$, $E/\Q$ be an elliptic curve and $\Delta_E$ be the minimal discriminant of $E$. Suppose that $\Q(E[m])=\Q(E[n])$. Then, for all primes $p$ such that $v_p(m)\neq v_p(n)$, we have $p\mid2\Delta_E$.
\end{corollaire}

\subsection{Index of images}\label{index of images}

Let $p$ be a prime. Let $M$ be a subgroup of $\mathrm{GL}_n(\Z_p)$ and let $G$ be a subgroup of $M$. In our setting we only need to consider the groups $\SL(\Z_p)$, $\GL(\Z_p)$ and $(\Z_p)^*$, but we state results in the general case as the approach is the same. For $k$ a positive integer, we denote by $M_k$ the image of $M$ in $\mathrm{GL}_n(\Z/p^k\Z)$ and $G_k$ the image of $G$ in $M_k$. We set $i_k=[M_k:G_{k}]$. We have $i_k\mid i_{k+1}$ by Remark~\ref{divisibilité des indices}. 
Moreover,

\begin{equation}\label{valeurs des suites d'index}\frac{i_{k+1}}{i_k}=\frac{\#M_{k+1}}{\# G_{k+1}}\cdot\frac{\# G_k}{\# M_k} \left| \frac{\# M_{k+1}}{\# M_k}\right.=\left\{\begin{array}{ll} p^{n^2} &\mbox{if } M=\mathrm{GL}_n(\Z_p)\\ p^{n^2-1}&\mbox{if } M=\mathrm{SL}_n(\Z_p).\end{array}\right.\end{equation}

from \eqref{noyau de la projection}. In particular, \begin{equation}\label{valeurs suites pour coincidence} G_k\simeq G_{k+1}\iff \frac{i_{k+1}}{i_k}=\left\{\begin{array}{ll} p^{n^2} &\mbox{if } M=\mathrm{GL}_n(\Z_p)\\ p^{n^2-1}&\mbox{if } M=\mathrm{SL}_n(\Z_p).\end{array}\right.\end{equation}

The idea of the lemma below and its proof follows~\cite[Lemma 3.7]{suth&zywina}.

\begin{prop}\label{suite i_k}Suppose that $M=\mathrm{GL}_n(\Z_p)$ or $\mathrm{SL}_n(\Z_p)$. 
The sequence $(u_k)=\left(\frac{i_{k+1}}{i_k}\right)$ satisfies $u_{k+1}\mid u_k$ for $k\geq1$ if $p$ is odd and for $k\geq2$ if $p=2$.
\end{prop}

\begin{proof}Suppose that $p$ is odd or that $k\geq2$ and $p=2$. Let $H_{k}$ be the kernel of the reduction map $G_{k}\to G_{k-1}$. Let $h\in G$ whose image in $G_{k}$ belongs to $H_{k}$. Then $h=I+p^{k-1}A$ with $A\in M_n(\Z_p)$. The map\[\begin{array}{c@{\;}c@{\;}c@{\;}c}
   \phi : & H_k & \longrightarrow & H_{k+1}\\
          & \overline{h} & \longmapsto & \overline{h^p} \end{array}\]
is an injective morphism since \[(I+p^{k-1}A)^p=I+\binom{p}{1}p^{k-1}A+\binom{p}{2}p^{2{k-1}}A^2+\dots\equiv I+p^{k}A\Mod{p^{k+1}}.\]Therefore $\frac{\#G_{k}}{\#G_{k-1}}$ divides $\frac{\# G_{k+1}}{\# G_{k}}$ and so, since $\frac{\# M_{k}}{\# M_{k-1}}=\frac{\# M_{k+1}}{\# M_{k}}$ from the equation~\eqref{valeurs des suites d'index}, we obtain $u_{k+1}\mid u_{k}$.
\end{proof}

\begin{corollaire}\label{split liftable from the beggining}
Let $k\geq1$ if $p$ is odd and $k\geq2$ if $p$ is even. If $M=\mathrm{GL}_n(\Z_p)$ or $M=\mathrm{SL}_n(\Z_p)$, and $G_k\simeq G_{k+1}$, then $G_1\simeq G_2\simeq\dots\simeq G_{k+1}$ if $p$ is odd, and $G_2\simeq G_3\simeq\dots\simeq G_{k+1}$ if $p$ is even.
\end{corollaire}

\begin{proof}
Suppose that $M=\mathrm{GL}_n(\Z_p)$. Equivalence \eqref{valeurs des suites d'index} gives $\frac{i_{k+1}}{i_k}=p^{n^2}$. Since the sequence $\left(\frac{i_{s+1}}{i_s}\right)$ is non-increasing from Lemma~\ref{suite i_k} for $s\geq1$ and $p$ odd or $s\geq2$ and $p=2$, and has values dividing $p^{n^2}$ by Equation~\eqref{valeurs suites pour coincidence}, then $\frac{i_{s+1}}{i_s}=p^4$ for all $s\leq k$. The proof is similar for $M=\mathrm{SL}_n(\Z_p)$.
\end{proof}

\begin{theorem}\label{indice coincidence det}Let $q=p$ and $k\geq1$ if $p$ is odd, or $q=p^2$ and $k\geq2$ if $p$ is even. Let $E/F$ be an elliptic curve. If $F(E[p^k])=F(E[p^{k+1}])$, then $F(E[q])=F(E[p^{k+1}])$.
\end{theorem}

\begin{proof}Let $M=\GL(\Z_p)$ and $G=\rho_{E,p^\infty}(G_F)$. So \[G_k=\rho_{E,p^k}(G_F)\simeq \Gal(F(E[p^k])/F).\]Therefore, the equality $F(E[p^k])=F(E[p^{k+1}])$ is equivalent to $G_k\simeq G_{k+1}$, and we use Corollary~\ref{split liftable from the beggining}.
\end{proof}

\begin{corollaire}\label{corollaire split et coincidence}
Let $k\geq2$ if $p$ is odd and $k\geq3$ if $p$ is even. Let $E/F$ be an elliptic curve such that the exact sequence\[\scalebox{0.95}{$1\to\Gal(F(E[p^{k+1}])/F(E[p^k]))\to\Gal(F(E[p^{k+1}])/F)\to\Gal(F(E[p^k])/F)\to1$}\]is split, then $F(E[p^{k-1}])=F(E[p^k])$.
\end{corollaire}

\begin{proof}
By Proposition~\ref{constructing coincidence from split}, there exists $L/F$ linearly disjoint from $F(E[p^k])/F$ such that~$L(E[p^k])=L(E[p^{k+1}])$. Therefore, by Theorem~\ref{indice coincidence det}, we have $L(E[p^{k-1}])=L(E[p^k])$. But $L\cap F(E[p^k])=F$ gives $L(E[p^{k-1}])\cap F(E[p^k])=F(E[p^{k-1}])$, from which we obtain \[[F(E[p^k]):F(E[p^{k-1}])]=[L(E[p^{k}]):L(E[p^{k-1}])]=1,\]using \cite[A.V.13, Proposition 5]{Bourbaki2}.
\end{proof}

\begin{corollaire}
If $p$ is odd or $k\geq3$, then the sequence \[\scalebox{0.95}{$1\to\Gal(\Q(E[p^{k+1}])/\Q(E[p^k]))\to\Gal(\Q(E[p^{k+1}])/\Q)\to\Gal(\Q(E[p^k])/\Q)\to1$}\]is not split.
\end{corollaire}

\begin{proof}
If the sequence is split, then the coincidence $\Q(E[p^k])=\Q(E[p^{k-1}])$ holds from Corollary~\ref{corollaire split et coincidence}, which is impossible by \cite[Theorem 1.4]{coincidences}.
\end{proof}

\begin{remark}
We notice that the assumptions of the theorem are necessary. Indeed, let $E/\Q$ be an elliptic curve satisfying $\Gal(\Q(E[8])/\Q)\simeq (\Z/2\Z)^6$. Then, the sequence in Corollary~\ref{corollaire split et coincidence} is split. However, since $\Gal(\Q(E[8])/\Q(E[4]))\leq(\Z/2\Z)^4$ by \eqref{noyau de la projection}, we have $(\Z/2\Z)^2\leq \Gal(\Q(E[4])/\Q)$ and so $\Q(E[2])\neq\Q(E[4])$.
\end{remark}

We obtain the following proposition on the adelic index:

\begin{prop}\label{p^4 divise l'indice}Let $E/F$ be an elliptic curve without CM, with a $(p^k,p^{k+1})$-coincidence. Then $[\GL(\hat{\Z}):\rho_E(G_F)]$ is divisible by $p^{4k}$ if $p$ is odd, or $\mathrm{max}\left\{2^4,2^{4k-1}\right\}$ if $p$ is even.
\end{prop}

\begin{proof}With the introduced notation, we consider $M=\GL(\Z_p)$ and $G=\rho_{E,p^\infty}(G_F)$. The index $i_{k+1}$ divides $[\GL(\hat{\Z}):\rho_E(G_F)]$ by Remark~\ref{divisibilité des indices} and so $i_{k+1}/i_k$ divides the global index $[\GL(\hat{\Z}):\rho_E(G_F)]$. But having a $(p^k,p^{k+1})$-coincidence is equivalent to have $G_k\simeq G_{k+1}$. Then, the proposition follows from Equation~\ref{valeurs suites pour coincidence} and Corollary~\ref{split liftable from the beggining}.
\end{proof}

\begin{remark}\label{indices Zywina}
In \cite[Theorem 1.3]{zywina2022possible}, Zywina gives a set of possible index $[\GL(\hat{\Z}):\rho_E(G_\Q)]$ which is generic for elliptic curves $E$ defined over $\Q$. In \cite[Theorem 1.2]{zywina2024numberfields}, he gives upper bounds for the index $[\GL(\hat{\Z}):\rho_E(G_F)]$ for elliptic curves $E/F$, depending on the number field $F$.
\end{remark}

\begin{remark}\label{coincidence et i_k}Let $E/F$ be an elliptic curve and take $G=\rho_{E,p^\infty}(G_F)$. So we have $G_k=\rho_{E,p^k}(G_F)$ for $k\geq1$.
\begin{enumerate}
    \item The sequence $(i_k)$ is increasing, and, if $E/F$ does not have CM, becomes stationary.

    \item For $E/\Q$ without CM, in \cite[Definition 2.21]{Zoe_yvon_these}, we define the $p$-adic depth $s$ as the smallest $k$ such that $\rho_{E,p^\infty}(G_{\Q})$ is the full inverse image of $\rho_{E,p^k}(G_\Q)$ in $\GL(p^\infty)$. Then, the sequence $\left(\frac{i_{k+1}}{i_k}\right)$ stabilizes at $1$ from $\mathrm{max}\{1,s\}$. In particular, if $G=\GL(\Z_p)$, then $\left(\frac{i_{k+1}}{i_k}\right)$ is constant equal to $1$. The converse is false: the elliptic curve with LMFDB label \href{https://www.lmfdb.org/EllipticCurve/Q/11/a/2}{11.a2} has non maximal Galois representation at $5$ and the sequence $\left(\frac{i_{k+1}}{i_k}\right)$ attached to $5$ is constant, equal to $1$. See \cite[Figure 2.1]{Zoe_yvon_these} for the possible $p$-adic depth for any $p$.
    
    \item  We consider $k\geq1$ if $p$ is odd, or $k\geq2$ if $p$ is even. If the first term of $\left(\frac{i_{k+1}}{i_k}\right)$ is different than $p^4$, then $E/F$ does not have $(p^k,p^{k+1})$-coincidences for any $k$. Otherwise, if $s$ is the rank of the first jump of the sequence, then $E/F$ does not have $(p^k,p^{k+1})$-coincidences for $k\geq s$.
    
    \item  We consider $k\geq1$ if $p$ is odd, or $k\geq2$ if $p$ is even. We know that the sequences $\left(\frac{i_{k+1}}{i_{k}}\right)$ is non-increasing, then constant. We can ask if it is decreasing then constant. The answer is no. For example, the elliptic curve with LMFDB label \href{https://www.lmfdb.org/EllipticCurve/Q/15/a/4}{15.a4} has, for $p=2$, $\left(\frac{i_{k+1}}{i_k}\right)=(2^2,2,2,2,1,\dots)$. We can also ask if the graphs are "progressively non-increasing", meaning that $\frac{i_{k+1}}{i_k}\in\left\{\frac{i_k}{i_{k-1}},\frac{1}{p}\frac{i_k}{i_{k-1}}\right\}$. The answer is also no. For example, the elliptic curve with LMFDB label \href{https://www.lmfdb.org/EllipticCurve/Q/15/a/8}{15.a8} has, for $p=2$, $\left(\frac{i_{k+1}}{i_k}\right)=(2^3,2,1,\dots)$ and the elliptic curve with LMFDB label \href{https://www.lmfdb.org/EllipticCurve/Q/40/a/4}{40.a4} has, for $p=2$, $\left(\frac{i_{k+1}}{i_k}\right)=(2^4,1,\dots)$.
\end{enumerate}
\end{remark}

\begin{example} 
Here are some examples over $\Q$, computed from \cite{lmfdb}, illustrating different possibilities for the sequence $\left(\frac{i_{k+1}}{i_k}\right)$. Using point (2) of Remark~\ref{coincidence et i_k} and \cite[Table 2.1]{Zoe_yvon_these}, the sequence $\left(\frac{i_{k+1}}{i_k}\right)$ is constant equal to~$1$ for $p\geq 13$, and it stabilizes at most at rank $5$ if $p=2$, at rank $3$ if $p=3$ and at rank $2$ if $p=5,7,11$.

\vspace*{0.5cm}
\centering{\scalebox{0.9}{\begin{tabular}{|c|c|c|l|}
  \hline 
  LMFDB & Minimal Weierstrass equation & Non max & \multicolumn{1}{c|}{Sequence $\left(\frac{i_{k+1}}{i_k}\right)$}\\
  label & & $p$ & \multicolumn{1}{c|}{attached to $p$}\\
  \hline
  \href{https://www.lmfdb.org/EllipticCurve/Q/14/a/6}{14.a6} & $y^2+xy+y=x^3+4x-6$ & $2$ & $1,2,1,\dots$\\
  \hline
  \href{https://www.lmfdb.org/EllipticCurve/Q/15/a/1}{15.a1} &$y^2+xy+y=x^3+x^2-2160x-39540$ & $2$ & $2^2,2^2,2,1,\dots$\\
  \hline
  \href{https://www.lmfdb.org/EllipticCurve/Q/15/a/2}{15.a2} & $y^2+xy+y=x^3+x^2-135x-660$ & $2$ & $2^2,2^2,1,\dots$\\
  \hline
  \href{https://www.lmfdb.org/EllipticCurve/Q/15/a/4}{15.a4} & $y^2+xy+y=x^3+x^2-80x+242$ & $2$ & $2^2,2,2,2,1,\dots$\\
  \hline
  \href{https://www.lmfdb.org/EllipticCurve/Q/15/a/5}{15.a5} & $y^2+xy+y=x^3+x^2-10x-10$ & $2$ & $2^3,2,1,\dots$\\
  \hline
  \href{https://www.lmfdb.org/EllipticCurve/Q/15/a/8}{15.a8} & $y^2+xy+y=x^3+x^2+35x-28$ & $2$ & $2^3,2^2,1,\dots$\\
  \hline
  \href{https://www.lmfdb.org/EllipticCurve/Q/20/a/3}{20.a3} & $y^2=x^3+x^2-x$ & $2$ & $2,2,1,\dots$\\
  \hline
  \href{https://www.lmfdb.org/EllipticCurve/Q/40/a/4}{40.a4} & $y^2=x^3+13x-34$ & $2$ & $2^4,1,\dots$\\
  \hline
  \href{https://www.lmfdb.org/EllipticCurve/Q/19/a/1}{19.a1} & $y^2+y=x^3+x^2-769x-8470$ & $3$ & $3,3,1,\dots$\\
  \hline 
  \href{https://www.lmfdb.org/EllipticCurve/Q/54/a/2}{54.a2} & $y^2+xy=x^3-x^2-3x+3$ & $3$ & $3^2,1,\dots$\\
  \hline
  \href{https://www.lmfdb.org/EllipticCurve/Q/11/a/1}{11.a1}& $y^2+y=x^3-x^2-7820x-263580$ & $5$ & $5,1,\dots$\\
  \hline
  \href{https://www.lmfdb.org/EllipticCurve/Q/11/a/2}{11.a2} & $y^2+y=x^3-x^2-10x-20$ & $5$ & $1,1,\dots$\\
  \hline
\end{tabular}}}

\justifying
\vspace*{0.5cm}
\end{example}

\begin{remark}\label{suite j_k}
Let $E/F$ be an elliptic curve and take $G=\rho_{E,p^\infty}(G_F)$. Then $\det G$ is a subgroup of $(\Z_p)^*$, with image $\det G_k$ in $(\Z/p^k\Z)^*$. We recall that $\det G_k\simeq\Gal(F(\zeta_{p^k})/F)$, by Proposition~\ref{weil pairing}. We set $j_k=[(\Z/p^k\Z)^*:\det G_k]$. We consider $k\geq1$ if $p$ is odd, and $k\geq2$ is $p$ is even. The sequence $\left(\frac{j_{k+1}}{j_k}\right)$ is non-increasing and has value in $\{1,p\}$ by Lemma~\ref{suite i_k} and the equation~\eqref{valeurs des suites d'index}. We have $j_{k+1}/j_k=p$ if and only if $\det G_k\simeq \det G_{k+1}$,  by the isomorphism~\eqref{valeurs suites pour coincidence}, and this is equivalent to $F(\zeta_{p^{k+1}})=F(\zeta_{p^k})$. Corollary~\ref{split liftable from the beggining} implies that $F(\zeta_{q})=F(\zeta_{p^{k+1}})$ with $q=p$ if $p$ is odd and $q=4$ if $p$ is even. As a consequence, the sequence $(j_k)$ is increasing and becomes stationary from the smallest $s$ such that $\zeta_{p^s}\notin F(\zeta_q)$.
\end{remark}

\begin{remark}\label{suite ell_k}Let $E/F$ be an elliptic curve and take $G=\rho_{E,p^\infty}(G_F)$. We observe that $\SL(\Z_p)\cap G$ is a subgroup of $\SL(Z_p)$ and so we can consider its projection in $\SL(\Z/p^k\Z)$ for each $k$. Unfortunately, these projections are not necessarily equal to $\SL(p^{k+1})\cap G_k$ and so we cannot use these groups to deal with the coincidence : if $G_k\simeq G_{k+1}$, we do not necessarily have $\SL(p^k)\cap G_k\simeq \SL(p^{k+1})\cap G_{k+1}$. Setting $\ell_k=[\SL(p^k):\SL(p^k)\cap G_k]$, we have $i_k=j_k\ell_k$ with $j_k$ defined as in Remark~\ref{suite j_k}. Suppose that $G_k\simeq G_{k+1}$. Then $i_{k+1}/i_k=p^4$. Hence \begin{align*}p^4=\frac{j_{k+1}\ell_{k+1}}{j_k\ell_k}&=p^4\frac{\#\det G_k}{\#\det G_{k+1}}\frac{\#(\SL(p^k)\cap G_k)}{\#(\SL(p^{k+1})\cap G_{k+1})}\\&=p^4\frac{[F(\zeta_{p^{k}}):F]}{[F(\zeta_{p^{k+1}}):F]}\frac{\#(\SL(p^k)\cap G_k)}{ \#(\SL(p^{k+1})\cap G_{k+1})}\cdot\end{align*}Then we have two situations:\begin{enumerate}
    \item $F(\zeta_{p^k})=F(\zeta_{p^{k+1}})$, and $\SL(p^k)\cap G_k\simeq \SL(p^{k+1})\cap G_{k+1}$,
    \item $F(\zeta_{p^{k+1}})\neq F(\zeta_{p^{k+1}})$, and the reduction map $\SL(p^{k+1})\cap G_{k+1}\to \SL(p^k)\cap G_k$ is not surjective.
\end{enumerate}
In the first case, we have seen in Remark~\ref{suite j_k} that $F(\zeta_{p^{k+1}})$ is equal to $F(\zeta_4)$ if $p=2$ and $F(\zeta_p)$ otherwise. The examples of vertical coincidence we have for elliptic curves over $\Q$ fit, obviously, in the second case. Indeed, the elliptic curve with LMFDB label \href{https://www.lmfdb.org/EllipticCurve/Q/40/a/4}{40.a4} has a $(2,4)$-coincidence with \[G_2=\rho_{E,4}(G_F)\simeq\left\langle\begin{pmatrix}0&1\\1&0\end{pmatrix}\right\rangle.\] In this case, we have $G_2\cap \SL(4)=\{\id\}$, whereas \[G_1\cap\SL(2)=G_1=\left\langle\begin{pmatrix}0&1\\1&0\end{pmatrix}\right\rangle.\]
\end{remark}

\subsection{Split liftable subgroups}\label{section split liftable}
Let $m,n$ be positive integers. If $E/F$ is an elliptic curve with an $(m,mn)$-coincidence, then $\rho_{E,m}(G_F)\simeq\rho_{E,mn}(G_F)$ and the image of $\rho_{E,mn}(G_F)$ in $\GL(m)$ is $\rho_{E,m}(G_F)$. It leads to the following definition. For a subgroup $G$ of $\GL(m)$, there are a priori several liftings of $G$ in $\GL(mn)$.

\begin{definition}
We say that a subgroup $G$ of $\GL(m)$ is \emph{split liftable} modulo $mn$ if there exists $G'\leq\GL(mn)$ such that $G$ is the image of $G'$ in $\GL(m)$ and $G\simeq G'$. We say that an element $g$ of $\GL(m)$ is \emph{split liftable} modulo $m$ if there exists $g'\in\GL(mn)$ with same order as $g$ and such that $g$ is the image of $g'$ in $\GL(m)$.
\end{definition}

A subgroup $G$ of $\GL(m)$ is split liftable modulo $mn$ is there exists an injective morphism $G\to\GL(mn)$ which makes the following diagram commutative:


\[\begin{tikzcd}
&\GL(mn)\arrow[d]\\
G \arrow[ur, hook] \arrow[r,hook] &\GL(m)
\end{tikzcd}\]
The definition above is up to conjugation, since two conjugate groups are isomorphic. Therefore:

\begin{prop}
Let $E/F$ be an elliptic curve with an $(m,mn)$-coincidence. Then $\rho_{E,m}(G_F)$ is split liftable modulo $mn$.
\end{prop}

The aim of this section is to determine the subgroups of $\GL(m)$ which are split liftable or not modulo some multiple of $m$.

\begin{remark}
If $G$ in $\GL(m)$ is split liftable modulo $mn$, then $G$ is split liftable modulo every $km$ such that $1\leq k\leq n$.
\end{remark}

\begin{remark}
In \cite{elkies3adic}, Elkies already use the property of being split liftable to construct the modular curve $\mathcal{X}_9$. It is defined by $\mathcal{X}_9=X(9)/(G/\left\langle-\id\right\rangle)$ where $G$ is a split lifting of $\SL(3)$ in $\GL(9)$, and more specifically, in $\SL(9)$.
\end{remark}

\begin{remark}
Corollary~\ref{split liftable from the beggining} tells that, if a subgroup of $\GL(p^k)$ is split liftable modulo $\GL(p^{k+1})$, then its image in $\GL(q)$ is also split liftable modulo $\GL(p^{k+1})$, where $q=p$ if $p$ is odd or $q=4$ if $p$ is even.
\end{remark}

\begin{lemma}\label{subgroups of split liftable subgroups}
Let $G\leq\GL(m)$ be split liftable modulo $mn$. Then, all subgroups of $G$ are split liftable modulo $mn$.
\end{lemma}

\begin{proof}Let $G'$ be a subgroup of $\GL(mn)$ such that $G$ is the image of $G'$ in $\GL(mn)$ and $G\simeq G'$. Then, the restriction $\pi:G'\to G$ of the natural projection $\GL(mn)\to\GL(m)$ is an isomorphism. Let $H\leq G$. We set $H':= \pi^{-1}(H)$. Then $H'\simeq H$ and $H$ is the image of $H'$ in $\GL(m)$.
\end{proof}

\begin{prop}\label{split liftable modulo every multiple}
Let $m\geq2$. The following subgroups of $\GL(m)$ are split liftable modulo every multiple of $m$:
\[\left\langle \begin{pmatrix}0&-1\\1&1\end{pmatrix},\begin{pmatrix}0&1\\1&0\end{pmatrix}\right\rangle,\quad\left\langle \begin{pmatrix}1&0\\0&-1\end{pmatrix}, \begin{pmatrix}0&1\\-1&0\end{pmatrix}\right\rangle.\]
\end{prop}

\begin{proof}Looking the first two groups as subgroups of $\GL(\Z)$, we observe that they have finite orders, respectively $12$ and $8$, and their elements only have coefficients in $\{0,1,-1\}$. Consequently, it is isomorphic to their projection modulo any integers $m$ such that $1\neq -1\Mod m$, that is for any $m\geq 3$. The case $m=2$ is given by Example~\ref{split liftable mod 2}. 
\end{proof}

\begin{example}\label{split liftable mod 2}The group $\GL(2)$ lifts to $\left\langle \begin{pmatrix}-1&1\\-1&0\end{pmatrix},\begin{pmatrix} 0&1\\1&0\end{pmatrix}\right\rangle\subseteq\GL(\Z)$, and so is split liftable modulo every even integer.
\end{example}

\begin{remark}
Let $E/F$ be an elliptic curve. To have an $(m,mn)$-coincidence, it is necessary to have $\rho_{E,m}(G_F)$ split liftable modulo $mn$, but it is not sufficient. Indeed, $\GL(2)$ is split liftable modulo $8$, but there are no $(2,8)$-coincidence for elliptic curve defined over $\Q$, see \cite[Theorem 1.4]{coincidences}. This is also the case for $\GL(3)$: the subgroup \[\left\langle\begin{pmatrix}1&0\\0&-1\end{pmatrix}, \begin{pmatrix}-2&2\\-2&-2\end{pmatrix},\begin{pmatrix}4&-2\\-3&4\end{pmatrix}\right\rangle\subseteq\GL(9)\]is a lifting of $\GL(3)$ of order $48$. Hence, $\GL(3)$ is split liftable modulo $9$ and yet there is no $(3,9)$-coincidence for elliptic curves with surjective mod $3$ Galois representation, by Corollary~\ref{vertical coincidence and trivial intersection} and Proposition~\ref{weil pairing}.
\end{remark}

Now we will present groups which are not split liftable, which give us obstructions having the coincidence. We underline that, by Lemma~\ref{subgroups of split liftable subgroups}, if $g\in \GL(m)$ is not split liftable modulo $mn$, then all groups containing $g$ are not split liftable modulo $mn$.

From now on, we will denote by $T$ the matrix $\begin{pmatrix}1&1\\0&1\end{pmatrix}$.



\begin{lemma}\label{borel subgroup}Let $p$ be a prime and $k\geq1$. The matrix $T$ in $\GL(p^k)$ is split liftable modulo $p^{k+1}$ if and only if $p=2,3$ and $k=1$.
\end{lemma}

\begin{proof}In $\GL(2)$, resp. $\GL(3)$, the matrix $T$ is conjugates to $\begin{pmatrix}0&1\\1&0\end{pmatrix}$, resp. $\begin{pmatrix}0&1\\-1&-1\end{pmatrix}$, and so is split liftable modulo $4$, resp. modulo $9$, by Lemma~\ref{split liftable modulo every multiple}. Set $q=p^2$ if $p=2,3$ and $q=p$ otherwise. We know, from Corollary~\ref{split liftable from the beggining}, that, for $k\geq1$ if $p$ is odd and $k\geq2$ if $p$ is even, if $T$ in $\GL(p^k)$ is split liftable modulo $p^{k+1}$, then $T$ in $\GL(q)$ is split liftable modulo $p^{k+1}$ and so modulo $pq$.
Now, if $T$ was split liftable modulo $pq$, then we could find $M\in M_2(\Z_p)$ such that $T+qM$ has order $q$ in $\GL(pq)$. But\[\scalebox{0.95}{$\left(\begin{smallmatrix}1+qa&1+qb\\qc&1+qd\end{smallmatrix}\right)^n\equiv\left(\begin{smallmatrix}1+nqa+\frac{n(n-1)}{2}qc\quad&n+nqb+\frac{n(n-1)}{2}q(a+d)+\frac{n(n-1)(n-2)}{6}qc\\nqc&1+nqd+\frac{n(n-1)}{2}qc\end{smallmatrix}\right)\Mod{pq}.$}\]Now, we take $n=q$. Then $p$ divides $\frac{q(q-1)}{2}$ and $\frac{q(q-1)(q-2)}{6}$. We obtain \[\left(T+qM\right)^q\equiv\begin{pmatrix}1&q\\0&1\end{pmatrix}\not\equiv \id\Mod{pq}. \]
\end{proof}

\begin{theorem}\label{T and vertical coincidence}
Let $p$ be a prime. Let $q=p$ if $p\neq 2,3$, or $q=p^2$ if $p=2,3$. Let $E/F$ be an elliptic curve. If $\rho_{E,q}(G_F)$ contains $T$, then $F(E[q])\neq F(E[pq])$.
\end{theorem}

In particular, if $E/F$ does not have CM and has a $(p^k,p^{k+1})$-coincidence, then $E/F$ has non maximal image modulo $p^2$, and even modulo $p$ if $p\neq 2,3$.

\begin{corollaire}
Let $E/F$ be an elliptic curve with multiplicative reduction at a prime $\got{r}$ of $\mathcal{O}_F$ and let $p$ be a prime not dividing $2v_\got{r}(j(E))$. Then $E/F$ does not have a $(p^k,p^{k+1})$-coincidence for any $k$.
\end{corollaire}

\begin{proof}
If $E/F$ has multiplicative reduction at a prime ideal $\got{r}$ and $p\nmid 2 v_\got{r}(j(E))$, then $T\in\rho_{E,p}(G_F)$ by \cite[Proposition 1.6]{ATAEC}. But this is not possible for a $(p^k,p^{k+1})$-coincidence from Theorem~\ref{T and vertical coincidence}.
\end{proof}

\begin{remark}Suppose that $p$ is odd, and $k\geq2$ if $p=3$. To study $(p^k,p^{k+1})$-coincidences further, it remains to deal with subgroups of $\GL(p^k)$ with non-surjective determinant, by Corollary~\ref{vertical coincidence and trivial intersection}, and which does not contains $T$, by Theorem~\ref{T and vertical coincidence}.
\end{remark}

\subsection{CM case}

If $E/F$ has complex multiplication by a quadratic field $K$ and $F\subseteq K(j(E))$, then we can say more.

\begin{prop}\label{coincidence CM}
Let $E/\overline{\Q}$ be an elliptic curve with CM by a quadratic field $K$ and $\K=K(j(E))$. If $\K(E[p^k])=\K(E[p^{k+1}])$, then $p=2$ and $k=1$.
\end{prop}

\begin{proof}
We have the following field inclusions:

\[\begin{tikzpicture}
\node      (A) {$\K(E[p^{k+1}])$};
\node   (R)       [below=of A] {$\quad$};
\node (B) [left=of R] {$\K(E[p^k])$};
\node (C) [right=of R]{$\K(h(E[p^{k+1}])$};
\node      (D)       [below=of R] {$\K(h(E[p^k]))$};

\draw[-] (A.south west) edge node[above] {$a$} (B.north);
\draw[-] (A.south east) edge node[above] {} (C.north);
\draw[-] (B.south) edge node[below] {$c$} (D.north west);
\draw[-] (C.south) edge node[below] {$d$}(D.north east);
\end{tikzpicture}\]
where $h$ is a Weber function for $E$ (see \cite{lozrobcm}).
Suppose that $k\geq1$ if $p$ is odd or $k\geq2$ if $p=2$ and $a=1$.
By \cite[Theorem 4.3]{lozrobcm}, we have, \[d=[\K(h(E[p^{k+1}])):\K(h(E[p^{k}]))]=p^2.\]
This implies that $p^2\mid c$. Moreover, by \cite[Theorem 4.1]{lozrobcm} we have $c\mid\# \mathcal{O}_K^*$. But $\# \mathcal{O}_K^*=2$, $4$ or $6$. Thus $p=2$ and $c=\# \mathcal{O}_K^*=4$. But $\mathcal{O}_K^*\simeq \Aut(E)$, so $j(E)=1728$ by \cite[III, Theorem 10.1]{AEC} and $\K=K(j(E))=K=\Q(i)$. In this case, $E$ is defined over $\Q$ and $\Q(E[2^k])\subsetneq \Q(E[2^{k+1}])$ by \cite[Proposition 3.9]{coincidences}. Moreover, the Weil pairing implies that\[\K=\Q(i)\subseteq \Q(E[4])\subseteq \Q(E[2^k])\subsetneq \Q(E[2^{k+1}])\] and so $\K(E[2^k])\subsetneq \K(E[2^{k+1}])$. We conclude that $p=2$ and $k=1$.
\end{proof}

\begin{remark}If an elliptic curve $E/F$ has a $(2,4)$-coincidence, then $\rho_{E,4}(G_F)$ must be a split lifting of $\rho_{E,2}(G_F)$. Example~\ref{split liftable mod 2} gives such a split lifting and a Magma computation shows that this is the only one up to conjugation (\cite[Proof of Proposition 3.9]{coincidences}). The corresponding modular curve is $X_{20b}$ in the notation of Rouse and Zureick-Brown \cite[Remark 1.6]{2adicimage}. They have computed its model, see \href{https://users.wfu.edu/rouseja/2adic/X20b.html}{https://users.wfu.edu/rouseja/2adic/X20b.html}, and so the map to the $j$-line, explicitely given in \cite[Proof of Proposition 3.9]{coincidences}. If the elliptic curve $E/F$ with $j$-invariant $j(E)$ has a $(2,4)$-coincidence, then there exists $t\in F$ such that \[j(E)=\frac{-4t^8+32t^7+80t^6-288t^5-504t^4+864t^3+1296t^2-864t-1188}{t^4+4t^3+6t^2+4t+1}\cdot\]For rational CM $j$-invariant, there is no such $t$.
\end{remark}

\section{Large images}\label{section 4}

Since $\GL(m)$ and $\GL(n)$ are not isomorphic for $m\neq n$, then $(m,n)$-coincidences cannot happen for elliptic curves with surjective mod $m$ and mod $n$ representations. In this section, we show that, under some conditions on $\K$, $(m,n)$-coincidences cannot happen if one of the images, said the image mod $m$, is large, \emph{i.e.} it contains the special linear group. In this case, the elliptic curve $E/F$ does not have CM, and it is said that it has maximal image at $m$. We will only deal with $m$ odd.

For a group $G$, we denote by $\D(G)$ its commutator subgroup. We know that $D(G)$ is the smallest normal subgroup of $G$ such that $G/\D(G)$ is abelian and this last is called the \emph{abelianization of $G$}. We will use several well-known results about the derived group of $\SL(m)$ and $\GL(m)$, given with detailed proofs in Appendix~\ref{appendix}.

We denote by $F^\mathrm{ab}$ the maximal abelian extension of $F$. We will compare the maximal abelian extension of $F(E[m])$ and that of $F(E[n])$.

\begin{prop}\label{abelianisé mod m}Let $m$ be an odd integer and $E/F$ be an elliptic curve. Suppose that $\rho_{E,m}(G_\K)$ contains $\SL(m)$. Then \[F(E[m])\cap F^\mathrm{ab}=\left\{\begin{array}{ll} F(\zeta_{m}) &\mbox{if } \mathrm{D}(\rho_{E,m}(G_\K))=\SL(m)\\ \text{a }\Z/3\Z\text{-extension of } F(\zeta_m)&\mbox{otherwise}.\end{array}\right.\]
\end{prop}

\begin{proof}Let $G=\rho_{E,m}(G_\K)$. We have $\SL(m)\leq G\leq \GL(m)$. Then, using Proposition~\ref{groupe dérivé GL}, \[\D(\SL(m))\leq \D(G)\leq \SL(m)).\]
Suppose that $\D(G)=\SL(m)$. Since \[G/\SL(m)\simeq\det(G),\] therefore the largest abelian quotient of $\K(E[m])/\K$ has Galois group isomorphic to $\det(G)$. By the Weil pairing, $\K(\zeta_{m})\subseteq\K(E[m])$ and from Proposition~\ref{weil pairing}, we have $\Gal(\K(\zeta_{m})/\K)\simeq\det(G)$. Then the largest abelian subextension of $\K(E[m])$ is $\K(\zeta_{m})$.

Now, suppose that $\D(G)\neq\SL(m)$. If $3\nmid m$, this does not happen, since in this case $\D(\SL(m))=\SL(m)$ from Proposition~\ref{groupe dérivé SL}. If $3\mid m$, then $\D(\SL(m))$ has index $3$ in $\SL(m)$ from Proposition~\ref{groupe dérivé SL} and so is $\D(G)$. It follows that $F(E[m])\cap F^\mathrm{ab}$ is an extension of degree $3$ of $F(\zeta_m)$.
\end{proof}

\begin{remark}\label{cas abélianisé plus grand}
The case $F(E[m])\cap F^\mathrm{ab}\neq F(\zeta_{m})$ occurs only for $\gcd(m,12)\neq 1$ by Proposition~\ref{groupe dérivé SL} and Proposition~\ref{groupe dérivé GL}. Let $k=v_3(m)$. In the previous proposition, $m$ is odd, and so $F(E[m])\cap F^\mathrm{ab}\neq F(\zeta_{m})$ only if $k>0$. In this case $L:=F(E[3^{k}])\cap F^\mathrm{ab}$ is a $(\Z/3\Z)$-extension of $F(\zeta_{3^{k}})$ and \[F(E[m])\cap F^\mathrm{ab}=L\otimes_F F(\zeta_{\frac{m}{3^k}}).\]
\end{remark}

\begin{remark}
Let $k=1$ if $p\geq5$, $k=2$ if $p=3$ and $k=3$ if $p=2$. If $\rho_{E,p^k}(G_F)$ contains $\SL(p)$, then $\rho_{E,p^{\infty}}(G_F)$ contains $\SL(\Z_p)$, by \cite[IV.3.4. Lemma~3]{SerreMG}.
\end{remark}

\begin{theorem}\label{large image m odd}
Let $m$ be an odd integer, $n\nmid m$ such that $\zeta_n\notin F$ and $E/F$ be an elliptic curve.
Suppose that $\rho_{E,m}(G_\K)$ contains $\SL(m)$ and that $E/F$ has an $(m,n)$-coincidence. Then \[3\mid m,\quad\text{and}\quad \mathrm{D}(\rho_{E,m}(G_F))\neq\SL(m),\quad \text{and}\quad F(\zeta_n)\subseteq L\] with $L$ a $\Z/3\Z$-extension of $F(\zeta_m)$.
\end{theorem}

\begin{proof}Suppose that $F(E[m])=F(E[n])$. Then $F(\zeta_n)\subseteq F(E[m])\cap F^\mathrm{ab}$. If $F(E[m])\cap F^\mathrm{ab}=F(\zeta_m)$, this is not possible since $\zeta_n\notin F$ and $n\nmid m$. Hence, from Proposition~\ref{abelianisé mod m}, $F(E[m])\cap F^\mathrm{ab}$ is a $\Z/3\Z$-extension of $F(\zeta_m)$ and the derived group of $\rho_{E,m}(G_F)$ is smaller than $\SL(m)$, which only happen if $3\mid m$ by Remark~\ref{cas abélianisé plus grand}.
\end{proof}

\begin{remark}
Under the hypotheses of the previous theorem, for $(m,n)=(p^k,p^{k+1})$, we know by the previous theorem that $p=2$ or $3$. Since $\SL(p^k)$ contains the matrix $T$, this result was already known by Theorem~\ref{T and vertical coincidence} and even more: the assumption $\zeta_{p^{k+1}}\notin F$ is unnecessary and $k=1$.
\end{remark}

By \cite[Lemma 1.7]{zywina2022explicit}, we have
$\SL(\hat{\Z})\cap\rho_E(G_F)=\D(\rho_E(G_F))$.
In particular, if $\rho_E(G_F)$ contains $\SL(\hat{\Z})$, then $\D(\rho_E(G_F))=\SL(\hat{\Z})$. In particular, if $8\cdot 9\mid m$ and $\SL(m)\subseteq \rho_{E,m}(G_F)$ then $\SL(\hat{\Z})\subseteq \rho_E(G_F)$ by \cite[Lemma 1 and the following paragraph]{antwerpIII}, and so $\D(\rho_{E,m})=\SL(m)$. Hence we deduce from Theorem~\ref{abelianisé mod m} the following:

\begin{theorem}\label{coincidence 72 divise m}
Let $E/F$ be an elliptic curve and $m$ be a positive integer. Suppose that $72\mid m$ and $\rho_{E,m}(G_F)$ contains $\SL(m)$. Then $F(E[m])\cap F^\mathrm{ab}=F(\zeta_m)$. In particular, if~$E/F$ has an $(m,n)$-coincidence, then $F(\zeta_n)\subseteq F(\zeta_m)$.
\end{theorem}

We finish the section by the following lemma, which contains some additional information about the extension $F(E[3])/F$.

\begin{lemma}\label{j(E) 3}
Let $E/F$ be an elliptic curve with $j$-invariant $j(E)$. Then we have $F(j(E)^{1/3})\subseteq F(E[3])$. Moreover, if $\SL(3)\leq\rho_{E,3}(G_F)$, then the extension $F(j(E)^{1/3})/F$ is non trivial.
\end{lemma}

\begin{proof}We know that we can parameterize the modular curve $X(1)$ with $j:X(1)\to\PP^1$ the $j$-invariant. Let $C_{ns}^+(3)$ be the normalizer of the non-split Cartan subgroup of $\GL(3)$ and the associate modular curve $X_{ns}^+(3)=X(3)/C_{ns}^+(3)$. From \cite[Proposition 4.1]{chen}, there exists a uniformizer $t:X_{ns}^+(3)\to\PP^1$ such that $t^3=j$. The points of $X_{ns}^+(3)(F(E[3]))$ correspond to elliptic curves defined over $F(E[3])$ with Galois image contained in $C_{ns}^+(3)$, that is every elliptic curves defined over $F$, since $\rho_{E,3}(G_{F(E[3])})=\id$. Hence $t(E)\in F(E[3])$. Using that $t^3=j$, we obtain $j(E)^{1/3}\in F(E[3])$. Moreover, if $j(E)^\frac{1}{3}\in F$, then $\rho_{E,3}(G_F)\subseteq C_{ns}^+(3)$. In particular, in this case, $\rho_{E,3}(G_F)$ cannot contain $\SL(3)$.
\end{proof}

\begin{remark}If $\zeta_3\in F$, then $F(j(E)^\frac{1}{3})/F$ is Galois and so is contained in $F^\mathrm{ab}$. If, moreover, $F(j(E)^\frac{1}{3})\cap F(\zeta_{m})=F$, then $F(E[m])\cap F^\mathrm{ab}=F(j(E)^\frac{1}{3},\zeta_m)$. Hence, if we have an $(m,n)$-coincidence, then \[3\mid m,\quad\text{and}\quad  \mathrm{D}(\rho_{E,m}(G_F))\neq\SL(m),\quad \text{and}\quad F(\zeta_n)\subseteq F(j(E)^\frac{1}{3},\zeta_m).\]
\end{remark}

\begin{remark}If $\zeta_3\in F$, $3\mid m$ and $\D(\rho_{E,3^{v_3(m)}}(G_\K))=\SL(3^{v_3(m)})$, then $F(j(E)^\frac{1}{3})\subseteq F(\zeta_{3^{v_3(m)}})$ and so we have $r\geq1$ such that \[F(j(E)^\frac{1}{3})=F(\zeta_{3^{r+1}})\quad \text{and}\quad F=F(\zeta_{3^r})\neq F(\zeta_{3^{r+1}}).\]
\end{remark}

\section{Coincidence of division fields of two elliptic curves}\label{section 5}

Let $E$ and $E'$ be elliptic curves defined over $F$. In all this section, we suppose that $F(E[m])=F(E'[n])$ for two integers $m$ and $n$. Then $\zeta_n\in F(E[m])$ and in the previous sections we gave constraints to this property. Hence:

\begin{prop}Let $p$ be a prime and $r$ be the largest integer such that $\Q(\zeta_{p^{r}})\subseteq \K\cap\Q(\mu_{p^\infty})$. Suppose that $p>q$ for all primes $q\mid m$ and $v_p(n)>r$. Then, $v_p(n)=1$ (and $r=0$), unless $(m,p)=(2^j,3)$ for some $j\geq1$, in which case $r=0$ and $v_p(n)\leq 2$, or $v_p(n)=r+1$.
\end{prop}

\begin{proof}
This follows from Proposition~\ref{greatest prime divisor} and the fact that $F(\zeta_{p^k})\subseteq F(E'[p^k])\subseteq F(E'[n])$.
\end{proof}

We can use Theorem~\ref{table of ramification division} and Remark~\ref{table of ramification cyclo} as the same way as in Section~\ref{section 2} to give constraints on a coincidence $F(E[m])=F(E'[n])$ using again that we must have $\zeta_n\in F(E[m])$. For example, we give a generalization of Corollary~\ref{coincidence and reduction}.

\begin{theorem}\label{coincidence and reduction two different ec}
For all primes $p$ such that $p\mid n$ and $p\nmid m\Delta_F$ we are in one of the following situation:
\begin{itemize}
    \item $v_p(n)=1$ and $E/F$ has bad reduction at every ideal above $p$,
    \item $v_p(n)=2$, $p=2$ and at each prime above $p$, $E/\K$ has either additive or non split multiplicative reduction,
    \item $v_p(n)=2$, $p=3$, and $E/\K$ has additive and potential good reduction at every ideal above $p$,
    \item $v_p(n)=3$ or $4$, $p=2$ and $E/\K$ has additive and potential good reduction at every ideal above $p$.
\end{itemize}
\end{theorem}

\begin{theorem}\label{vertical coincidence and trivial intersection two ec}
Let $p$ be a prime and $k\geq1$ such that $F\cap\Q(\zeta_{p^k})=\Q$. If $\K(E[p^k])=\K(E'[p^{k+1}])$, then $p=2$.
\end{theorem}

\begin{proof}
This follows from Theorem~\ref{equality for p odd}.
\end{proof}

\begin{theorem}For all primes $p$ such that $v_p(m)\neq v_p(n)$, we have \[p\mid 2\cdot\Delta_F\cdot\mathrm{N}(\got{f}_E)\cdot \mathrm{N}(\got{f}_{E'}).\]
\end{theorem}

\begin{proof}
The proof is the same as for Theorem~\ref{coincidence ramified or bad red}, using Proposition~\ref{coincidence and reduction two different ec} and Theorem~\ref{vertical coincidence and trivial intersection two ec} instead of Corollary~\ref{coincidence and reduction} and Corollary~\ref{vertical coincidence and trivial intersection}.
\end{proof}

\begin{theorem}
Suppose that $m$ is odd, $n\nmid m$, $\zeta_n\notin F$ and $\rho_{E,m}(G_\K)$ contains $\SL(m)$. Then \[3\mid m,\quad\text{and}\quad \mathrm{D}(\rho_{E,m}(G_F))\neq\SL(m),\quad \text{and}\quad F(\zeta_n)\subseteq L\] with $L$ a $\Z/3\Z$-extension of $F(\zeta_m)$.
\end{theorem}

\begin{proof}
The proof is exactly the same as Theorem~\ref{large image m odd}, replacing $F(E[n])$ by $F(E'[n])$.
\end{proof}

\begin{remark}
Except in Subsection~\ref{section 3.2}, the results of Section~\ref{section 3} use the reduction link between $\rho_{E,p^{k+1}}$ and $\rho_{E,p^k}$ on the same elliptic curve. In particular, this method does not apply for different elliptic curves.
\end{remark}


\bigskip

\bibliographystyle{annotate}
\bibliography{Ref}

\appendix

\section{Derived groups of $\GL(m)$ and $\SL(m)$}\label{appendix}

In this appendix, we give elementary and detailed proofs of well-known results about the derived groups of $\GL(m)$ and $\SL(m)$, for any integer $m$. They are used in Section~\ref{section 4} in the case where $m$ is odd. For a group $G$, we denote by $\D(G)$ its commutator subgroup, generated by all the elements $[g,h]=ghg^{-1}h^{-1}$ with $g,h\in G$. We know that $\D(G)$ is normal in $G$ and is the smallest group such that $G/\D(G)$ is an abelian group: the \emph{abelianization of $G$}.

We recall that $\SL(\Z)$ is generated by $S=\begin{pmatrix}0&-1\\1&0\end{pmatrix}$ and $T=\begin{pmatrix}1&1\\0&1\end{pmatrix}$. They satisfy $S^2=(ST)^3=-I$. 




\begin{prop}\label{gen gr dérivé Z}
The quotient group $\SL(\Z)/\D(\SL(\Z))$ is cyclic of order $12$, generated by the equivalence class of $T$.
\end{prop}

\begin{proof}
Let $\bar{S}$ and $\bar{T}$ the classes of $S$ and $T$ in $\SL(\Z)/\mathrm{D}(\SL(\Z))$ respectively. Since $\SL(\Z)/\mathrm{D}(\SL(\Z))$ is abelian, we have $\bar{S}\bar{T}=\bar{T}\bar{S}$. Hence $(\bar{S}\bar{T})^3=\bar{S}^3\bar{T}^3=\bar{S}^2$, which gives $\bar{S}=\bar{T}^{-3}$ and $\bar{T}^{12}=\bar{S}^4=I$. It follows that $\bar{T}$ generated $\SL(\Z)/\mathrm{D}(\SL(\Z))$ and has order $12$.
\end{proof}

Let $m\mid n$ be positive integers and, for $i=n,m$, set $X_i=\SL(i))/\mathrm{D}(\SL(i))$.
The following diagram has exact rows and is commutative.
    \begin{equation}\label{diagram gr dérivé}\xymatrix{
    1 \ar[r] & \D(\SL(\Z)) \ar[r] \ar[d] & \SL(\Z) \ar[r]\ar@{->>}[d] & \Z/12\Z \ar[r]\ar[d] & 0 \\
    1 \ar[r] & \D(\SL(n)) \ar[r]\ar[d] & \SL(n) \ar[r]\ar@{->>}[d] & X_n \ar[r]\ar[d] & 0 \\
     1 \ar[r] & \D(\SL(m)) \ar[r] & \SL(m) \ar[r] & X_m \ar[r] & 0 
  }\end{equation}

\begin{prop}\label{groupe dérivé SL}
Let $m$ be a positive integer. The abelianization of $\SL(m)$ is isomorphic to $\Z/\gcd(m,12)\Z$. In particular, if $m$ is coprime to $6$, then $\SL(m)$ is perfect.
\end{prop}

\begin{proof}By Diagram~\eqref{diagram gr dérivé} and Proposition~\ref{gen gr dérivé Z}, the image of $T$ in $\SL(m)$ generates the abelianization of $\SL(m)$, whose order divides $12$. Moreover, the image of $T$ in $\SL(m)$ has order $m$. Hence the abelianization of $\SL(m)$ has order dividing $\gcd(m,12)$. To prove that its order is exactly $\gcd(m,12)$, it suffices to prove it for $m=2,3$ and $4$, and then to use again Diagram~\eqref{diagram gr dérivé}. We define the commutators \[A=\left[\begin{pmatrix}-1&-1\\0&-1\end{pmatrix},\begin{pmatrix}1&-2\\1&-1\end{pmatrix}\right]=\begin{pmatrix}3&-1\\1&0\end{pmatrix},\]
\[B=\left[\begin{pmatrix}-1&-1\\0&-1\end{pmatrix},\begin{pmatrix}0&-1\\1&0\end{pmatrix}\right]=\begin{pmatrix}2&1\\1&1\end{pmatrix}.\]
Let $G=\langle A,B\rangle\subseteq\SL(\Z)$. For $m=2,3,4$, let $G_m$ be the image of $G$ is $\SL(m)$. By computing explicitely $G_2$, $G_3$ and $G_4$, we find that they have respectively order $3$, $8$ and $12$, and they are normal subgroups of $\SL(m)$. Then $\SL(m)/G_m$ is an abelian group of order $m$. Hence, as a subgroup of $\Z/12\Z$, it is isomorphic to $\Z/m\Z=\Z/\gcd(m,12)\Z$.
\end{proof}

\begin{remark}
We also can prove the previous result for $m=p^k$ with $p\geq5$ prime, just by observing that the diagram ~\ref{diagram gr dérivé} gives that the order of $X_m$ is both a divisor of $12$ and a divisor of
$\#\SL(p^k)=p^{3(k-1)+1}(p-1)(p+1)$. Hence, the order of the abelianization of $\SL(p^k)$ is $1$, and so its derived group is itself, unless $p=2$ or $3$.
\end{remark}

\begin{prop}\label{groupe dérivé GL}
For an odd integer $m$, the derived group of $\GL(m)$ is $\SL(m)$, and so its abelianization is isomorphic to $(\Z/m\Z)^*$. If $m$ is even, the derived group of $\GL(m)$ has index $2$ in $\SL(m)$.
\end{prop}

\begin{proof}
For $A,B\in\GL(m)$, we have $\det(ABA^{-1}B^{-1})=1$. So $[A,B]\in\SL(m)$. Therefore \[\D(\SL(m))\leq\D(\GL(m))\leq\SL(m).\]
For $m$ coprime to $6$, we have proven that $\D(\SL(m))=\SL(m)$, and so $\D(\GL(m))=\SL(m)$. For $m=3^k$, we know that $\mathrm{D}(\SL(3^k))$ has index $3$ in $\SL(3^k)$, and so $\mathrm{D}(\GL(3^k))$ is either $\SL(3^k)$ or $\mathrm{D}(\SL(3^k))$. For $k=1$, $\mathrm{D}(\SL(3))$ is explicitly know by the proof of Proposition~\ref{groupe dérivé SL} and\[\left[\begin{pmatrix}1&-1\\1&1\end{pmatrix},\begin{pmatrix}1&1\\0&-1\end{pmatrix}\right]=\begin{pmatrix}1&1\\-1&0\end{pmatrix}\notin\mathrm{D}(\SL(3)).\]It follows that $\mathrm{D}(\GL(3))=\SL(3)$. We obtain the diagram below:

\[\begin{tikzcd}[column sep=1.2cm]
\D(\SL(3^k)) \arrow[r, hook]\arrow[d, two heads] & \D(\GL(3^k)) \arrow[r, hook]\arrow[d, two heads] & \SL(3^k) \arrow[d, two heads] \\
    \D(\SL(3)) \arrow[hook]{r}[below]{\text{index }3}  & \D(\GL(3)) \arrow[r, equal] & \SL(3)
\end{tikzcd}\]
showing that $\D(\SL(3^k))\neq \D(\GL(3^k))$ and so $\D(\GL(3^k))=\SL(3^k)$. For $m=2^k$, we use the same strategy. From Proposition~\ref{groupe dérivé SL}, we already know that $\D(\GL(2))$ has index $2$ in $\SL(2)$. We have the following diagram:
\[\begin{tikzcd}[column sep=1.2cm]
    \D(\SL(2^k)) \arrow[r, hook]\arrow[d, two heads] & \D(\GL(2^k)) \arrow[r, hook]\arrow[d, two heads] & \SL(2^k) \arrow[d, two heads] \\
    \D(\SL(4)) \arrow[r, hook]\arrow[d, two heads] & \D(\GL(4)) \arrow[r, hook]\arrow[d, two heads] & \SL(4) \arrow[d, two heads] \\
    \D(\SL(2)) \arrow[r, equal] & \D(\GL(2)) \arrow[hook]{r}[below]{\text{index }2} & \SL(2)
\end{tikzcd}\]
The group $\D(\SL(4))$ is known by the proof of Proposition~\ref{groupe dérivé SL}, and we have \[\left[\begin{pmatrix}0&1\\1&0\end{pmatrix},\begin{pmatrix}2&1\\-1&1\end{pmatrix}\right]=\begin{pmatrix}0&-1\\1&1\end{pmatrix}\in \D(\GL(4))\backslash \D(\SL(4)).\] 
Using that $\D(\SL(2^k))$ has index $4$ in $\SL(2^k)$ for $k\geq2$, we obtain that $\D(\GL(4))$ has exactly index $2$ in $\SL(4)$ and $\D(\GL(2^k))$ has exactly index $2$ in $\SL(2^k)$ for all $k\geq2$. The result follows, since the Chinese remainder theorem gives \[\D(\GL(m))\simeq\underset{\underset{p\text{ prime}}{p^k\mid\mid m}}{\prod} \D(\GL(p^k)).\]
\end{proof}



\end{document}